\documentclass[preprint,11pt]{elsarticle}
\usepackage[latin1]{inputenc}
\usepackage[english]{babel}
\usepackage{amssymb,amsmath,graphicx}
\usepackage{multirow}
\usepackage{natbib}
\usepackage{enumitem}
\usepackage{setspace}
\usepackage[table]{xcolor}
\usepackage{microtype}
\usepackage{varwidth}
\usepackage[noend]{algorithmic}
\usepackage{algorithm}
\usepackage[hidelinks]{hyperref}
\usepackage[all]{nowidow}
\usepackage{longtable}

\newcommand{\cG}{{\cal G}}
\newcommand{\cV}{{\cal V}}

\newcommand{\cE}{{\cal E}}

\newcolumntype{P}[1]{>{\centering}p{#1}}
\newcolumntype{M}[1]{>{\centering}m{#1}}
\newcolumntype{H}{>{\setbox0=\hbox\bgroup}c<{\egroup}@{}}
\definecolor{light-gray}{gray}{0.95}

\newcommand{\myblue}[1]{#1}

\allowdisplaybreaks

\def\sqw{\hbox{\rlap{\leavevmode\raise.3ex\hbox{$\sqcap$}}$%
\sqcup$}}
\def\sqb{\hbox{\hskip5pt\vrule width4pt height6pt depth1.5pt%
\hskip1pt}}

\def\qed{\ifmmode\hbox{\hfill\sqb}\else{\ifhmode\unskip\fi%
\nobreak\hfil
\penalty50\hskip1em\null\nobreak\hfil\sqb
\parfillskip=0pt\finalhyphendemerits=0\endgraf}\fi}
\def\cqfd{\ifmmode\sqw\else{\ifhmode\unskip\fi\nobreak\hfil
\penalty50\hskip1em\null\nobreak\hfil\sqw
\parfillskip=0pt\finalhyphendemerits=0\endgraf}\fi}

\title{The Vehicle Routing Problem with Service Level Constraints}

\author{Teobaldo Bulh\~oes, Minh Ho\`ang H\`a, Thibaut Vidal}

\begin{document}

\begin{center}

\vspace*{-0.5cm}

\begin{huge}
The Vehicle Routing Problem with \vspace*{0.2cm} \linebreak Service Level Constraints
\end{huge}

\vspace*{0.45cm}

\textbf{Teobaldo Bulh\~oes *} \\
Instituto de Computa\c c\~ao, Universidade Federal Fluminense, Niter\'oi, Brazil \\ 
tbulhoes@ic.uff.br\\
\vspace*{0.15cm}
\textbf{Minh Ho\`ang H\`a} \\
FPT Technology Research Institute, FPT University, 
Hanoi, Viet Nam \\
hoanghm@fpt.edu.vn\\
\vspace*{0.15cm}
\myblue{\textbf{Rafael Martinelli}} \\
\myblue{Departamento de Engenharia Industrial, Pontifícia~Universidade~Católica~do~Rio~de~Janeiro,~Brazil} \\
\myblue{martinelli@puc-rio.br} \\
\vspace*{0.15cm}
\textbf{Thibaut Vidal} \\
Departamento de Informática, Pontifícia Universidade Católica do Rio de Janeiro, Brazil \\
vidalt@inf.puc-rio.br \\
\vspace*{0.15cm}

\end{center}
\noindent
\textbf{Abstract.}
We consider a vehicle routing problem which seeks \myblue{to minimize cost} subject to service level constraints on several groups of deliveries. This problem captures some essential challenges faced by a logistics provider which operates transportation services for a limited number of partners and should respect contractual obligations on service levels. The problem also generalizes several important classes of vehicle routing problems with profits.
\myblue{To solve it, we propose a compact mathematical formulation, a branch-and-price algorithm, and a hybrid genetic algorithm with population management, which relies on} problem-tailored solution representation, crossover and local search operators, as well as an adaptive penalization mechanism \myblue{establishing} a good balance between service levels and costs. Our computational experiments
 show that the proposed heuristic returns very high-quality solutions for this difficult problem, matches all optimal solutions found for small and medium-scale benchmark instances, and improves upon existing algorithms for two important special cases: the vehicle routing problem with private fleet and common carrier, and the capacitated profitable tour problem. \myblue{The branch-and-price algorithm also produces new optimal solutions for all three problems.}
\vspace*{0.2cm}


\section{Introduction}
\label{section:intro}

In a recent industrial application, the authors have been confronted with a complex variant of \myblue{the} vehicle routing problem (VRP) which received, until now, only limited attention in the academic literature.
We take the viewpoint of a third-party logistics provider (\myblue{3PL}), which operates long-haul transportation \myblue{services} for a number of business partners. The \myblue{company} operates on a planning horizon and delivers products to various delivery locations as requested by the partners a few days in advance. A strict delivery deadline, in the form of a last possible delivery day, is set for each transportation request. The \myblue{company} has established agreements with each partner specifying, among others, a minimum level of on-time deliveries for its group of requests.

\myblue{The efficient management of collaborative logistics has stimulated a rich set of studies over the years \citep{jayaram-2010,leuschner-2014,stefansson-2006}. In our case,} 
even in the presence of precise information on delivery requests, the \myblue{company} faces the following optimization challenge: \myblue{because of limited available resources, i.e., fleet size and vehicle capacity,} all services cannot be realistically fulfilled. It is necessary to select a subset of deliveries and determine cost-efficient vehicle routes in such a way that the overall activity is profitable and that the agreements with the partners are respected. A natural strategy of the company consists in attempting to fulfill the service levels on a rolling horizon with some safety margin, hence ensuring overall satisfaction of contractual clauses on a larger time period (e.g., one month of activity). As such, the \myblue{firm} seeks to balance quality of service and operational costs, subject to a minimum threshold on quality of service for some groups of requests. \myblue{These group requirements create linking constraints between the service selection decisions, which need to be carefully considered during routing optimization.}

We now introduce a simple and deterministic variant of the VRP which captures some essential decisions in this situation. The study of this simplified problem will help to identify key properties and methods.
The VRP with service levels (VRP-SL) can be formulated as follows.
Let $\cG=(\cV, \cE)$ be a complete undirected graph with $|\cV|=n+1$ nodes. 
The node $v_0 \in \cV$ represents a depot, where a fleet of $m$ identical vehicles is based.
Each other node $v_i$ for $i \in \{1, \ldots , n\}$ represents a customer, associated with a demand~$q_i$, a profit~$p_i$, and a service weight~$s_i$ which represents its \myblue{relative importance} in the group service level constraint.

The set of customers is distributed into $K$ subsets: $\cV-\{v_0\} = \bigcup_{k=1,\dots,K} \cV_k$, such that $\cV_k \cap \cV_{k'} = \varnothing$ for any $k \neq k'$.  Each subset represents the deliveries of one partner and is associated with a requested service level $\alpha_k$.
Any edge $(i,j) \in \cE$ represents a possible trip between a node $v_i \in \cV$ and a node $v_j \in \cV$ \myblue{with} a distance cost $d_{ij}$. The goal of the VRP-SL is to find \emph{up to}~$m$ vehicle routes starting and ending at the depot, such that
\begin{itemize}[nosep]
\item[--] each customer is serviced at most one time,
\item[--] the total demand quantity of any route does not exceed a vehicle capacity $Q$,
\item[--] \myblue{the service level of each group $k$ is attained, i.e., the total service weight of the deliveries to this group reaches $\alpha_k \sum_{v_i \in \cV_k} s_i$, and}
\item[--] the sum of travel costs and lost profits is minimized.
\end{itemize}

This problem belongs to the wide class of vehicle routing problems with profits, which also includes the team orienteering problem (TOP), the profitable VRP (VRPP), the VRP with private fleet and common carrier (VRPPFCC) and the capacitated profitable tour problem (CPTP). Interestingly, as highlighted in Section \ref{section:Litt}, this problem fills a gap in the literature, since most known multi-vehicle problems with customer selections either aim to maximize service levels subject to distance constraints (TOP) or seek a weighted optimization of distance and service levels, through penalties for outsourcing or lost profits (VRPPFCC and CPTP). To this date, very few works on deterministic settings \citep{Tang2006a,Yadollahpour2009,Li2016} have addressed multi-vehicle routing optimization subject to a service level (SL) constraint.
Finally, the VRP-SL is finally a natural extension of the generalized VRP (GVRP, see \cite{Ghiani2000,Baldacci2009a,Bektas2011,Ha2014}), which also models a rich set of VRP applications.

To address the VRP-SL, we introduce a compact ILP formulation \myblue{and a branch-and-price algorithm} which can solve to optimality small- and medium-scale instances, as well as a hybrid population metaheuristic inspired by the \myblue{Unified Hybrid Genetic Search (UHGS)} framework of \cite{Vidal2012,Vidal2012b}. Previous applications of UHGS have led to efficient algorithms for several VRPs, including the GVRP, TOP, VRPP and VRPPFCC, among others~\citep{Vidal2012b,Vidal2014}. \myblue{However}, UHGS \myblue{relies heavily on} the fact that problem objectives and constraints are either cumulated on or separately applied to each route, hence allowing to optimize customer-selection decisions independently via a dedicated route evaluation operator. This decomposition does not apply to the VRP-SL without a Lagrangian relaxation of the group constraints, which would make it impossible to find the optimal solution in some cases and impede the overall method performance. Because of these characteristics, we had to revise most of the operators of the method while keeping the general principles. We thus use a two-chromosome solution representation, which contains a service level chromosome and a \myblue{giant-tour} chromosome. We derive a new crossover, use dedicated local search moves as well as a penalization strategy to represent, optimize and inherit customer-selection decisions. The contributions of this article are the following:
\begin{enumerate}[nosep,leftmargin=0.6cm]
\item We introduce the VRP-SL, a rich VRP connected with important applications in collaborative logistics, generalizing many classes of routing problems with profits.
\item We propose a compact mathematical formulation for the problem, a \myblue{branch-and-price algorithm}, as well as an efficient hybrid genetic search (HGS) which exploits problem-tailored selection representation, crossover, local searches, and penalty management operators to find a good balance between service levels and costs.
\item We conduct extensive computational experiments to evaluate the performance of the proposed methods on new benchmark instances for the VRP-SL as well as classical benchmark instances for the VRPPFCC and CPTP. The proposed HGS finds all known optimal solutions for the considered problems, outperforms previous methods for the VRPPFCC and CPTP, and \myblue{generates solutions of consistent quality on large instances. Several solutions of the VRPPFCC and CPTP are also proven optimal for the first time by the branch-and-price algorithm.}
\end{enumerate}

The remaining parts of the article are organized as follows.
\myblue{
Section~\ref{section:Litt} reviews the related literature.
Section~\ref{section:math} presents some mathematical formulations for the problem. Sections \ref{section:BP} and \ref{section:HGA} describe the branch-and-price and the hybrid genetic algorithm, respectively. Section~\ref{section:experiments} reports the experimental analyses, and Section~\ref{section:concl} concludes.
}

\section{Related literature and subproblems}
\label{section:Litt}

\begin{table}[!htpb]
\renewcommand{\arraystretch}{1.15}
\vspace*{-0.5cm}
\hspace*{0.3cm}
 \rotatebox{90}{%
 \begin{varwidth}{1.1\textheight}
\scalebox{0.82}
{
\begin{tabular}{|l|m{7.5cm}@{\hspace*{0.1cm}}r|m{7.5cm}@{\hspace*{0.1cm}}r|m{7.5cm}@{\hspace*{0.1cm}}r|}
\cline{2-7}
\multicolumn{1}{c|}{\strut}\vspace*{-0.35cm}&&&&&&\\
\multicolumn{1}{c|}{\strut}&\multicolumn{2}{c|}{\textbf{Min $\alpha \times$cost $- \ \beta \times$service level}}&\multicolumn{2}{c|}{\textbf{Max service level s.t. cost constraints}}&\multicolumn{2}{c|}{\textbf{Min cost s.t. service level constraints}} \\
\multicolumn{1}{c|}{\strut}\vspace*{-0.35cm}&&&&&&\\
\hline
\vspace*{-0.25cm}&&&&&&\\
&\multicolumn{2}{c|}{\textbf{Ib) Profitable tour problem (PTP)}}&\multicolumn{2}{c|}{\textbf{IIb) Orienteering problem (OP)}}&\multicolumn{2}{c|}{\textbf{IIIb) Prize-collecting TSP (PCTSP)}}\\
\vspace*{-0.2cm}&&&&&&\\
&\emph{No dedicated heuristics, but:}&&\textbf{Heuristic:}&&\textbf{Heuristic:}&\\
&\multicolumn{2}{l|}{Reductions to elementary \myblue{shortest path} with }&GRASP&\cite{Campos2014}&Tabu search using GENIUS neighborhood&\cite{Pedro2013}\\ 
& possible negative cycles, or asymmetric TSP &&Guided local search&\cite{Vansteenwegen2009b}&Clustering search&\cite{Chaves2008}\\
\textbf{T}&Various approximation algorithms & \cite{Archetti2014}&Ant colony optimization and VNS&\cite{Schilde2009}&Lagrangian heuristic&\cite{DellAmico1995}\\
\textbf{S}&Extensions into PCTSP with profits&&Tabu search&\cite{Gendreau1998b}&Tabu Search (steel industry)&\cite{Lopez1998}\\
\textbf{P}& considered in the objective&&and others...&&&\\
&&&&&&\\
&\textbf{Exact:}&&\textbf{Exact:}&&\textbf{Exact:}&\\
&Branch-and-cut (+ capacity constraint)&\cite{Jepsen2014}&Branch-and-cut&\cite{Fischetti1998}&Branch-and-cut&\cite{Berube2009}\\
&Bounding method via Lagrangian relax.&\cite{DellAmico1995}&Branch-and-cut + Lagrangian relaxation&\cite{Kataoka1988}&Bounding method via Lagrangian relax.&\cite{DellAmico1995}\\
&&&&&Formulation + Additive bounds&\cite{Fischetti1988}\\
\vspace*{-0.25cm}&&&&&&\\
\hline
\vspace*{-0.25cm}&&&&&&\\
&\multicolumn{2}{c|}{\textbf{Ia) VRPPFCC and CPTP}}&\multicolumn{2}{c|}{\textbf{IIa) Team-orienteering problem (TOP)}}&\multicolumn{2}{c|}{\cellcolor{light-gray} \textbf{IIIa) Prize collecting VRP (PCVRP)}} \\
\vspace*{-0.2cm}&&&&&\cellcolor{light-gray}&\cellcolor{light-gray}\\
&\textbf{Heuristic VRPPFCC:}&&\textbf{Heuristic:}&&\multicolumn{2}{c|}{\cellcolor{light-gray}\emph{Limited number of works on this variant}} \\
&Hybrid GA with large neighborhoods&\cite{Vidal2014}&Pareto mimic algorithm&\cite{Ke2015}&\cellcolor{light-gray}&\cellcolor{light-gray} \\ 
&Adaptive variable neighborhood search&\cite{Huijink2014}&Hybrid GA with large neighborhoods&\cite{Vidal2014}&\cellcolor{light-gray}\textbf{Heuristic:}&\cellcolor{light-gray} \\ 
&Adaptive variable neighborhood search&\cite{Stenger2012a}&Multi-start simulated annealing&\cite{Lin2013}& \cellcolor{light-gray}Adaptive VNS for a related problem&\cellcolor{light-gray} \\
\textbf{V}&Tabu with possible ejection chains&\cite{Cote2009a,Potvin2011}&Augmented large neighborhood search&\cite{Kim2013}&\cellcolor{light-gray}with a minimum delivery level&\cellcolor{light-gray}\hspace*{0.3cm}\cite{Stenger2013}\\
\textbf{R}&and others...&&Particle swarm optimization&\cite{Dang2013}&\cellcolor{light-gray}Guided local search (steel industry)&\cellcolor{light-gray}\cite{Yadollahpour2009} \\
\textbf{P}&&&Path relinking&\cite{Souffriau2010}&\cellcolor{light-gray}Iterated local search (steel industry)&\cellcolor{light-gray}\cite{Tang2006a} \\
&\textbf{Heuristic CPTP:}&& Guided local search&\cite{Vansteenwegen2009b}&\cellcolor{light-gray} Two-stage self-adaptive VNS&\cellcolor{light-gray}\cite{Li2016}\\
&VNS and tabu search&\cite{Archetti2008d}&Ant colony optimization&\cite{Ke2008}&\cellcolor{light-gray}&\cellcolor{light-gray}\\
&&&and others...&&\cellcolor{light-gray}&\cellcolor{light-gray}\\
&&&&&\cellcolor{light-gray}&\cellcolor{light-gray}\\
&\textbf{Exact:}&&\textbf{Exact:}&&\cellcolor{light-gray}\textbf{Exact:}&\cellcolor{light-gray}\\
&Branch-and-price&\cite{Archetti2013a,Archetti2008d}&Branch-and-price&\cite{Archetti2013a,Archetti2008d,Boussier2007}&\cellcolor{light-gray}\emph{None to our knowledge}&\cellcolor{light-gray}\\
&&&Branch-and-cut-and-price&\cite{Keshtkaran2015}&\cellcolor{light-gray}&\cellcolor{light-gray} \\
&&&\myblue{Branch-and-cut}&\myblue{\cite{ElHajj2016}}&\cellcolor{light-gray}&\cellcolor{light-gray} \\
\vspace*{-0.25cm}&&&&&&\\
\hline
\end{tabular}
}
 \caption{Synthesis of state-of-the-art heuristics and exact approaches for single- and multi-vehicle routing problems with profits.}\label{Table:Summary}
 \end{varwidth}}
\end{table}

Vehicle routing problems with profits are the \myblue{subject} of an extensive literature, surveyed in \cite{Archetti2014,Feillet2005,Vansteenwegen2010,Vidal2012a}. Generally, one seeks to jointly minimize cost and maximize customer's service levels, two objectives which are conflicting when the profits do not render all deliveries profitable, or in the presence of additional vehicle constraints (e.g., distance or capacity limits). The related literature can thus be classified according to three fundamental solution techniques for this bi-objective~problem:
\begin{itemize}[nosep]
\item[\textbf{I)}] \textbf{Weighted sum --} minimizing a weighted difference of costs \mbox{and service levels};
\item[\textbf{II)}] \textbf{Constraints on cost/distance --} maximizing service levels subject to distance constraints (independently for each route);
\item[\textbf{III)}] \textbf{Constraints on service levels --} minimizing distance subject to service level constraints.
\end{itemize}

Table \ref{Table:Summary} presents an overview of these main classes of methods, for the single and multi-vehicle routing problems with profits (variants of the TSP and VRP, respectively).
Objectives I~and~II are the \myblue{subjects} of a vast literature, which includes a very large variety of exact methods and metaheuristics. The TOP, in particular, seeks the maximization of service levels subject to distance constraints, and has been studied in dozens of articles in the past decade (see \citep{Archetti2013a,Archetti2008d,Boussier2007,Dang2013,ElHajj2016,Ke2008,Ke2015,Kim2013,Lin2013,Souffriau2010,Vansteenwegen2009b,Vidal2014}, among others).
In contrast, a larger methodological gap remains for objective III, which aims to minimize travel cost subject to service level constraints. The multi-vehicle case of this objective has been mostly studied in the context of hot strip mill scheduling for steel production \cite{Yadollahpour2009,Zhang2009}. More recently, \cite{Li2016} investigated a special case with one group of customers (all of them) and one service level constraint, usually \myblue{referred to as the} prize-collecting VRP (PCVRP). The authors proposed a self-adaptive VNS and reported computational experiments on new instances. These instances, however, are now unavailable. The VRP-SL generalizes this problem by introducing \myblue{multiple} groups and associating profits to deliveries. By doing so, the problem remains concise and simple to define, but generalizes many VRP classes:
\begin{itemize}[nosep,leftmargin=0.6cm]
\item[--] The PCTSP and PCVRP correspond to VRP-SL instances with a single group.
\item[--] The CPTP, PTP and VRPPFCC all correspond to instances with 0\% service levels.
\item[--] The CVRP corresponds to instances with 100\% service levels.
\item[--] The generalized VRP (GVRP) can be reduced to the VRP-SL by imposing a small service level for each group, hence forcing at least one delivery.
\item[--] Any instance of the periodic vehicle routing problem (PVRP) such that, for each customer~$i$, any combination of $f_i$ different days in a set $D_i$ of available days is feasible, can be transformed into a VRP-SL instance with $\sum_i |D_i|$ customers. This is done by duplicating each~customer~$i$ into as many vertices as possible visit days, defining a group for the resulting visits with a service level $\alpha_i = f_i / |D_i|$, and setting a large cost $M$ for the edges that link two vertices from different days.
This reduction follows the same principles as the reduction from PVRP instances with frequency $f_i = 1$ into GVRP instances, discussed in \cite{Baldacci2009a}, but is more general as it allows to deal with frequency values greater than 1.
\end{itemize}

Finally, since the name \emph{prize-collecting VRP} is not consistently used in the literature and far from self-explanatory, we opted for the name \emph{VRP with service levels} (VRP-SL) for the proposed problem with several groups, which is clearer and more distinctive.

\section{Mathematical Formulations}
\label{section:math}

\noindent
\textbf{Compact formulation.}
We first propose a mixed integer linear formulation of the \mbox{VRP-SL}. Our mathematical model is based on the two-commodity flow formulation of \citep{Baldacci2004}, which has already been extended to several closely related VRP variants in  \cite{Ha2013, Ha2014}. The graph $\cG$ is first extended into $\overline{\cG}=(\overline{\cV}, \overline{\cE})$ by adding a new vertex~$v_{n+1}$, representing a copy of the depot. We thus define $\overline{\cV}=\cV \cup \{v_{n+1}\}$, $\cV^\prime = \overline{\cV}\setminus \{v_0, v_{n+1}\}$, $\overline{\cE}=\cE \cup \{(v_i, v_{n+1}),v_i \in \cV^\prime \}$, and $d_{i,n+1}=d_{0,i}$ for all $v_i \in \cV^\prime$.
\myblue{
For each edge $(v_i,v_j) \in \overline{\cE}$, we define a binary variable $x_{ij}$, set to $1$ if and only if a vehicle travels on this edge, as well as two flow variables $f_{ij}$ and $f_{ji}$. 
When a vehicle travels from $v_i$ to~$v_j$, the flow $f_{ij}$ represents the current load in the vehicle,
and the flow $f_{ji}$ represents the residual capacity of the vehicle ($f_{ji} = Q - f_{ij}$). Finally, $y_i$ is a binary variable which is set to $1$ if and only if $v_i \in \cV \setminus \{v_0\}$ is serviced.
The VRP-SL can be formulated as:}
\begin{align}
 \qquad \textrm{Minimize} \hspace*{1.5cm} \sum_{(v_i, v_j) \in \overline{\cE}}d_{ij}x_{ij} & + \sum_{v_i \in \cV \setminus \{v_0\}} p_i(1-y_i) \label{eq:obj1}\\
\textrm{Subject to} \hspace*{2.1cm} \sum_{v_i \in \cV_k}s_iy_i &\geq \alpha_k \sum_{v_i \in \cV_k}s_i &  k \in \{1,\dots,K\} \label{eq:c1}\\
\sum_{v_i \in \overline{\cV},i<k}x_{ik}+\sum_{v_j \in \overline{\cV},j>k}x_{kj}&=2y_k &  v_k \in \cV^\prime \label{eq:c2}\\
\sum_{v_j \in \overline{\cV}}(f_{ji}-f_{ij})&=2q_iy_i &  v_i \in \cV^\prime \label{eq:c3}\\
\sum_{v_j \in \cV^\prime}f_{0j}&=\sum_{v_i \in \cV^\prime}q_iy_i \label{eq:c4}\\
\sum_{v_j \in \cV^\prime}f_{n+1j}&=zQ \label{eq:c5}\\
f_{ij}+f_{ji}&=Qx_{ij} &  (v_i,v_j) \in \overline{\cE} \label{eq:c6}\\
x_{ij}& \in \{0, 1\} &  (v_i,v_j) \in \cE \label{eq:c8}\\
y_i & \in \{0, 1\} &  v_i \in \cV\setminus \{v_0\} \label{eq:c9}\\
f_{ij} &\geq 0,f_{ji} \geq 0 &  (v_i,v_j) \in \overline{\cE}\label{eq:c10}\\
z &\leq m, z \in \mathbb{N}  \label{eq:c11}
\end{align}

The objective of Equation (\ref{eq:obj1}) aims to minimize the sum of transportation costs and lost profits.
\myblue{
With this objective, the optimal solution value is always non-negative (no negative terms). Moreover, the optimal CVRP solution is a feasible solution of this model, with a cost greater or equal to the VRP-SL optimum.}
Constraints (\ref{eq:c1}) impose the service levels for each group. Constraints (\ref{eq:c2}) ensure that each vertex of $\cV \setminus \{v_0\}$ is visited \myblue{at most once. These constraints also connect the path variables~($x_i$) to the customer selection variables ($y_i$) in order to evaluate the profits}. \myblue{Constraints (\ref{eq:c3})--(\ref{eq:c6}) define a feasible two-commodity flow between the source $v_0$ and the sink~$v_{n+1}$.}
Specifically, Constraints~(\ref{eq:c3}) state that the inflow minus the outflow at each vertex $v_i \in \cV^\prime$ is equal to \myblue{$2 q_i$} if~$v_i$ is used, and $0$ otherwise.
The outflow at the source vertex~$v_0$, computed in Constraint~(\ref{eq:c4}), is set to the total demand of the vertices that are serviced in the solution,
and the outflow at the sink~$v_{n+1}$, calculated in Constraint~(\ref{eq:c5}), corresponds to the total capacity of the vehicle fleet.
\myblue{Finally, Constraints~(\ref{eq:c6}) establish the link between the flows $f_{ij}$ and $f_{ji}$, and Constraint~(\ref{eq:c11}) sets a bound on the number of vehicles.}

The linear relaxation of the VRP-SL can be strengthened via some simple valid inequalities. If $s_i$ is integer for all $v_i \in \cV \setminus \{v_0\}$, then the service level constraints (\ref{eq:c1}) can be transformed into:
\begin{align}
\sum_{v_i \in \cV_k}s_iy_i &\geq \Big\lceil\alpha_k \sum_{v_i \in \cV_k}s_i\Big\rceil & k \in \{1,\dots,K\}. \label{eq:c1eq}
\end{align}
Moreover, the following flow inequalities from \citep[Equation 64]{Baldacci2004} are used: 
\begin{align}
f_{ij} \geq q_j x_{ij} \text{ and } &f_{ji} \geq q_ix_{ij} & i, j \neq v_0 \text{ and } i, j \neq v_{n+1}. \label{eq:flow}
\end{align}
\myblue{The first inequality explicitly forces $f_{ij}$ to contain the delivery quantity $q_j$ for the next customer, and the second inequality forces the residual capacity $f_{ji}$ to be greater than the delivery quantity~$q_i$ of the previous customer.}
\myblue{
Finally, for each group $k \in \{1,\dots,K\}$, we define $Z_k$ as the minimum load quantity which allows to satisfy the service level constraint:
\begin{equation}
Z_k = \min_{y_i  \in \{0, 1\}}  \left\{  \sum_{v_i \in \cV_k} q_iy_i \hspace*{0.2cm} \bigg| \hspace*{0.2cm}   \sum_{v_i \in \cV_k}s_iy_i \geq \alpha_k \sum_{v_i \in \cV_k}s_i  \right \}
\end{equation}
These values can be evaluated in pseudo-polynomial time. Then, the following capacity cut is valid for each group $k$:
\begin{align}
\sum_{v_i \in V_k, v_j \in \cV \setminus \cV_k}x_{ij} &\geq \myblue{2} \left\lceil\frac{Z_k}{Q}\right\rceil & k \in \{1,\dots,K\}. \label{eq:capa}
\end{align}

Note that we also investigated in preliminary experiments a similar formulation based on one-commodity flows \citep{Gavish1978} for directed graphs, hence associating two variables $x_{ij}$ and~$x_{ji}$ for each edge ($v_i, v_j$). Although both formulations should produce similar bounds \citep{Letchford2006a}, the suggested two-commodity flow formulation of the VRP-SL led to overall better results in our context, possibly due to the smaller number of binary variables.\\}

\myblue{
\noindent
\textbf{Set Partitioning Formulation.}
We now present a set-partitioning based formulation of the VRP-SL, which corresponds to a Dantzig-Wolfe decomposition of the model of Equations (\ref{eq:obj1})--(\ref{eq:c11}). Similar formulations have been successfully used in the VRP literature in the past years to obtain stronger linear relaxations \citep{baldacci-2008,contardo-2014,fukasawa-2006,martinelli-2011}. The drawback of this formulation comes from its exponential number of variables, which must be tackled using a column generation algorithm.
}

\myblue{
Let $\Omega$ be the set of all feasible routes for the problem. A route $r \in \Omega$ is a closed walk starting and ending at the depot, visiting a set of customers only once  and respecting the vehicle's capacity (this definition will be extended afterwards to allow visiting a customer more than once). Then for each route, we define a binary variable $\lambda_r$ indicating whether the route $r$ is used in the solution. The resulting formulation is as follows:
\begin{align}
 \qquad \textrm{Minimize} \hspace*{1.5cm} \sum_{r \in \Omega}c_r \lambda_r & + \sum_{v_i \in \cV^\prime} p_i(1-y_i) \label{eq:sp-obj}\\
\textrm{Subject to} \hspace*{1.4cm} \sum_{v_i \in \cV_k}s_iy_i &\geq \alpha_k \sum_{v_i \in \cV_k}s_i &  k \in \{1,\dots,K\} \label{eq:sp-c1}\\
\sum_{r \in \Omega}a^r_i \lambda_r &= y_i & v_i \in \cV^\prime \label{eq:sp-c2}\\
\sum_{r \in \Omega}\lambda_r &\leq m & \label{eq:sp-c3}\\
\lambda_r & \in \{0, 1\} & r \in \Omega \label{eq:sp-c4}\\
y_i & \in \{0, 1\} &  v_i \in \cV^\prime \label{eq:sp-c5}
\end{align}

The objective function \eqref{eq:sp-obj} minimizes the total cost of the active routes plus the lost profits. Constraints \eqref{eq:sp-c1} are the same as Constraints \eqref{eq:c1}. In Constraints \eqref{eq:sp-c2}, $a^r_i$ (boolean) represents the number of times route $r$ visits customer $i$. Then each constraint forces one of the routes which visits customer $i$ to be active if variable $y_i$ is set to $1$. Constraints \eqref{eq:sp-c3} limit the number of active routes to the number of available vehicles.
}

\myblue{
\section{Branch-and-price algorithm}
\label{section:BP}
}

\myblue{
The set partitioning formulation has an exponential number of variables. To circumvent this issue, we use a column generation algorithm which starts with no variables and solves a pricing sub-problem to generate new ones. For the formulation presented in the last section, the pricing sub-problem is an Elementary Shortest Path Problem with Resource Constraints (ESPPRC). Ideally, one would want to price only elementary routes, but since the ESPPRC is known to be strongly NP-hard \cite{dror-1994}, we solve a pseudo-polynomial relaxation, the Shortest Path Problem with Resource Constraints (SPPRC) \cite{christofides-1981}, allowing the routes to visit a given customer more than once. For both problems, the objective is to find a route with negative reduced~cost. Given the dual variables $\beta_i$ and $\gamma$ associated to Constraints \eqref{eq:sp-c2} and \eqref{eq:sp-c3}, a constant $b^r_{ij}$ representing the number of times route $r$ traverses edge $(v_i, v_j)$ and setting $\beta_0 = \gamma$, the original reduced cost of a route and its reformulation as a function of the edges are stated in Equations~\eqref{eq:red-cost-1}~and~\eqref{eq:red-cost-2}:
\begin{align}
\overline{c}_r &= c_r - \gamma - \sum_{i \in \cV^\prime} a^r_i \beta_i
\label{eq:red-cost-1} \\
\overline{c}_r &= \sum_{(v_i, v_j) \in \cE} b^r_{ij}\left(c_{ij} - \frac{\beta_i + \beta_j}{2} \right). \label{eq:red-cost-2}
\end{align}

When relaxing the ESPPRC into the SPPRC, the bounds obtained by  column generation deteriorate in most cases. Some alternative relaxations have thus been proposed to find a better balance between efficiency and solution quality \cite{christofides-1981,irnich-2006,martinelli-2014}. We use the {\em ng}-route relaxation~\cite{baldacci-2011}, which has been successfully applied to multiple VRP variants in the last years. For each vertex~$i$, a set $NG_i \subseteq \cV^\prime$ is defined to represent the ``memory'' of $i$. These $NG_i$ sets usually contain a subset of vertices closest from $i$. The pricing sub-problem is solved by a forward dynamic programing algorithm, which maintains a memory of past visits to prohibit some vertices, but uses the $NG_i$ sets to reduce this memory size and thus the size of the state space. During the algorithm, each path $P$ has an associated label $\mathcal{L}(P)$ containing the last customer visited $v(P)$, the total reduced cost $\overline{c}(P)$, the current load $q(P)$ and the customers which have been visited and \emph{remembered} $\Pi(P)$. When extending a path from customer $v(P)$ to customer $v_j$, the extension is only allowed if $j \notin \Pi(P)$, and the label for the new path $P^\prime$ can be obtained as:}
\myblue{
\begin{equation}
\mathcal{L}(P^\prime) = \left(v_j, \overline{c}(P) + \overline{c}_{ij}, q(P) + q_{v_j}, \Pi(P) \cap NG_j \cup \{v_j\} \right)
\label{eq:path}
\end{equation}
Lastly, we use a dominance rule to fathom labels which cannot lead to an optimal solution. A label $\mathcal{L}(P_1)$ dominates a label $\mathcal{L}(P_2)$ if 
$$
\{v(P_1) = v(P_2)\} \wedge \{\overline{c}(P_1) \leq \overline{c}(P_2)\}  \wedge \{q(P_1) \leq q(P_2)\}  \wedge  \{\Pi(P_1) \subseteq \Pi(P_2)\}.
$$

Other techniques can be used to further improve the overall column generation efficiency. We use two approaches. The first one is a simple heuristic pricing. This algorithm stores only the label with best reduced cost for each customer $v_i$ and load $q$ during the dynamic programming. The second approach is dual stabilization \cite{pessoa-2010}. We use a parameter $\alpha \in [0, 1[$, as shown in Equations \eqref{eq:dual-stab-g}--\eqref{eq:dual-stab-b}, to avoid a large variation in the values of the dual variables between two iterations of the column generation. Our algorithm starts with $\alpha = 0.9$ and reduces the value of $\alpha$ by $0.1$ each time the pricing sub-problem returns an invalid route.
\begin{align}
\gamma &= \alpha \gamma^{k - 1} + (1 - \alpha) \gamma^k \label{eq:dual-stab-g}\\
\beta_i &= \alpha \beta^{k - 1}_i + (1 - \alpha) \beta^k_i & \forall v_i \in \cV^\prime \label{eq:dual-stab-b}
\end{align}

Finally, to improve the bounds obtained by the column generation algorithm, we embedded it into
a branch-and-bound (B\&B) procedure. At each node of the branch-and-bound tree, the column generation algorithm is called to obtain the solution of the linear relaxation, taking into account possible fixed variables due to branching. This algorithm framework is usually called branch-and-price (B\&P). 
A key difference with classic B\&B algorithms relates to the fact that it is not possible, in our context, to branch on the $\lambda_r$ variables, since a fixing of $\lambda_r = 0$ would result in repricing the variable. To overcome this difficulty, the B\&P branches on the original $x_{ij}$ variables as well as the $y_i$ variables, choosing at each node on the tree the most fractional variable to branch, but always giving priority to $y_i$ variables. Fixing a $y_i$ variable is straightforward. The $x_{ij}$ variables, however, are not explicitly present in the formulation, and thus we add Equation \eqref{eq:branch} to the set partitioning formulation for each fixed edge. In this equation, $b^r_{ij}$ represents the number of times route $r$ traverses edge $(v_i, v_j)$, and $\overline{x}_{ij}$ is the value which should be fixed for $x_{ij}$.
\begin{equation}
\sum_{r \in \Omega} b^r_{ij} \lambda_{r} = \overline{x}_{ij}
\label{eq:branch}
\end{equation}
The presence of Equation \eqref{eq:branch} introduces a new dual variable $\rho_{ij}$ which must be considered when computing the reduced costs, thus changing Equation \eqref{eq:red-cost-2} into Equation \eqref{eq:red-cost-bp}.
\begin{equation}
\overline{c}_r = \sum_{(v_i, v_j) \in \cE} b^r_{ij}\left(c_{ij} - \frac{\beta_i + \beta_j}{2} - \rho_{ij} \right)
\label{eq:red-cost-bp}
\end{equation}

This combination of techniques leads to an efficient exact method, whose performance will be analyzed in Section \ref{section:experiments}. \vspace*{0.3cm}}

\section{Population-based metaheuristic}
\label{section:HGA}

\myblue{The VRP-SL is known to be NP-hard as it generalizes the CVRP, and even sophisticated exact approaches can only solve small- and medium-scale problem instances within a reasonable CPU time. To fill this methodological gap and to solve the larger instances which arise in practice, we introduce a dedicated hybrid genetic search with advanced diversity control.}
 
\subsection{General structure of the method}

The method, illustrated in Algorithm \ref{alg:hgsadc}, uses the same resolution strategy as the unified hybrid genetic search (UHGS) of \cite{Vidal2012,Vidal2012b}.
Starting from an initial population, it iteratively selects two parents to generate an offspring individual via a crossover operator. This offspring is improved by means of a local search procedure and inserted in the population. This sequence of operations is performed until the termination of the method, once $\mathit{It}_{\text{NI}}$ successive iterations without improvement have been performed.

\begin{algorithm}[htbp]
\begin{spacing}{1.1}
\begin{algorithmic}[1]
\STATE \strut Initialize the population with random solutions
\WHILE {not $\mathit{It}_{\textsc{ni}}$ consecutive iterations without improvement of the best solution}
  \STATE Select two parents $P_1$ and $P_2$
  \STATE Generate an offspring $C$ by applying the crossover on $P_1$ and $P_2$
  \STATE Educate $C$ using local search
  \STATE Insert $C$ into the population
  \IF {$C$ is not feasible} 
  \STATE With 50\% probability, repair $C$ and insert it into the population
  \ENDIF
  \IF {\myblue{$\mathit{It}_{\textsc{div}}$ iterations elapsed since the last diversification}}
  \STATE Diversify the population
  \ENDIF
\ENDWHILE
\STATE \textbf{return} best feasible solution
\end{algorithmic}
\caption{\strut Hybrid Genetic Search  (HGS) for the VRP--SL}	  \label{alg:hgsadc}
\end{spacing}
\end{algorithm}

The method exploits penalized infeasible solutions in the population, which is divided into two sub-populations of feasible and infeasible individuals. Whenever one sub-population reaches a maximum size, a number of individuals are \myblue{eliminated to retain the best solutions}. A repair procedure is also applied on infeasible solutions to restore feasibility and generate additional feasible individuals. Finally, the approach uses an adaptive diversity management, which has been shown to be particularly successful when solving VRPs~\citep{Vidal2012}. Parents and survivors selections are driven by two criteria, cost and diversity contribution, rather than \myblue{solely on cost as in} traditional GAs, and additional diversification phases are implemented after \myblue{every} $\mathit{It}_{\textsc{div}}$ iterations without improvement to provide new solution characteristics to the population.

Finally, the proposed algorithm also significantly differs from UHGS in the definition of its basic building blocks: solution representation, crossover, local search moves, distance measure between individuals, and penalties allowed. These components are described below.

\subsection{Solution Representation and Evaluation}
\label{section:sol_repre}

\begin{figure}[htbp]
\centering
 \includegraphics[width=0.72\textwidth]{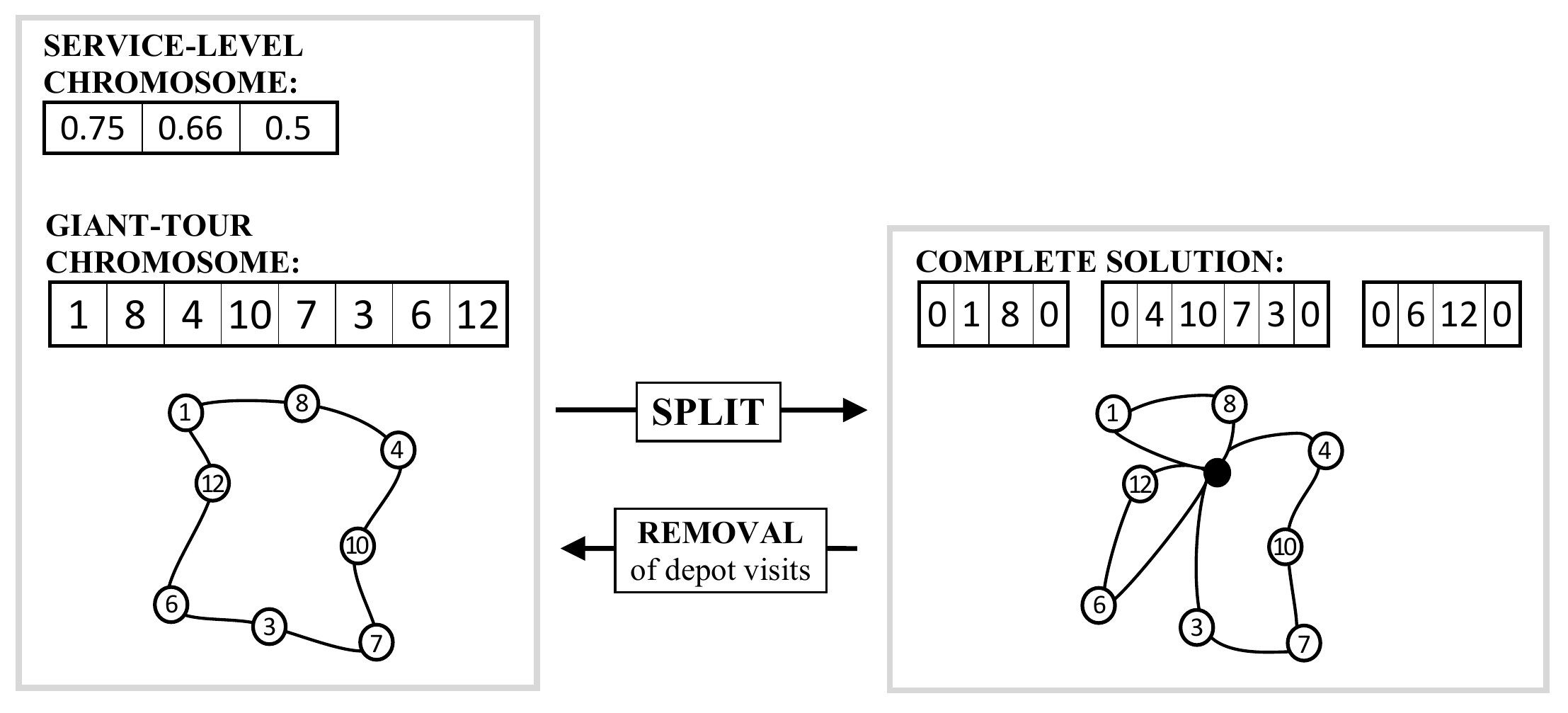}
 \caption{Solution representation and decoding via the Split algorithm. To compute the service levels, three groups are considered in the example: $\mathcal{V}_1 = \{1,2,3,4\}$, $\mathcal{V}_2 = \{5,6,7,8,9,10\}$ and $\mathcal{V}_3 = \{11,12\}$. For all $i$, $s_i = 1$.}
 \label{representation}
\end{figure}

In the proposed HGS, each individual in the population is represented by two chromosomes: a \emph{service level chromosome}, which gives the current service level of each group $k \in \{1,\dots,K\}$, and the \emph{\myblue{giant-tour} chromosome}, which provides a permutation of visits for the serviced customers, \myblue{without occurrences of the depot}.
\myblue{This solution representation is illustrated in Figure \ref{representation}. It is incomplete in the sense that some additional information, the locations of the visits to the depot, is needed to perform cost evaluations. Still, this information can be quickly recovered by means of a dynamic-programming algorithm, called  \emph{Split}, which optimally subdivides the \myblue{giant-tour} chromosome into separate routes. In the proposed HGS, the Split algorithm strictly respects the sequence of visits, i.e., it cannot select or exclude customer visits and modify the service levels. This is the classical context of application of the algorithm of \cite{Beasley1983} and \cite{Prins2004}, which reduces the splitting problem into the search of a shortest path in an acyclic graph in which each arc represents a possible route, i.e., a sequence of consecutive visits in the giant tour, connected to the depot. The reader is referred to \cite{Prins2004} and \cite{Vidal2016} for a detailed description of the Split algorithm as well as an efficient $\mathcal{O}(n)$ implementation.}

We now describe the cost function used for route and solution evaluations.
As observed in our computational experiments and in \cite{Vidal2013a}, the exploration of penalized infeasible solutions during the search has a significantly positive impact on solution quality. Thus, the cost $\phi(r)$ of a route $r$ involves the distance, the total profit associated to the customers which are visited, and a possible excess of load in a route multiplied by a penalty factor $w^Q$. Let $\phi^D(r) = \sum_{i=1}^{|r|-1} d_{r(i),r(i+1)}$, $\phi^Q(r) = \sum_{i=1}^{|r|} q_{r(i)}$ and $\phi^P(r) = \sum_{i=1}^{|r|} p_{r(i)}$ be, respectively, the total distance, load and profit collected in route $r$, then the route cost is defined as:
\begin{equation}
 \phi(r) = \phi^D(r) - \phi^P(r) + w^Q\max \{ 0, \phi^Q(r)-Q \}.
\end{equation}

We also explore infeasible solutions with respect to service level constraints. These linking constraints involve the whole solution, and thus the penalty is defined at the level of the solution evaluation rather than the route cost. \myblue{The penalized cost $\phi^\textsc{cost}(S)$ of a solution $S$, described as a set of routes $r \in S$, can thus be evaluated as:}
\begin{align}
 \phi^\textsc{cost}(S) &= \phi_P + \sum_{r \in S} \phi(r) + \sum_{k=1}^{K} w^S_k \max(\alpha_k - \phi^k(S), 0), \label{solution-cost} \\
\text{with } \phi^k(S) &= \sum_{v_i \in \cV_k \cap S} s_i \bigg/ \sum_{v_i \in \cV_k} s_i & k \in \{1,\dots,K\}.
\end{align}
In this equation, $\phi_P = \sum_{i=1}^n p_i$ is the total profit of all customers (a constant),
 $\phi^k(S)$ is the weight ratio of group $k$ in solution $S$, and $w^S_k$ for $k \in \{1,\ldots, K\}$ are the penalty coefficients associated to service level violations. These penalty coefficients are automatically adjusted by the method during the search, as explained in Section~\ref{section:pop}.

\subsection{Generation of New Individuals}
\label{section:generation}

Each new solution is generated by a successive application of the \emph{Selection}, \emph{Crossover}, and \emph{Education} operators, followed by a possible \emph{Repair}.

\begin{figure}[!ht]
\centering
\vspace*{0.1cm}
 \includegraphics[width=0.63\textwidth]{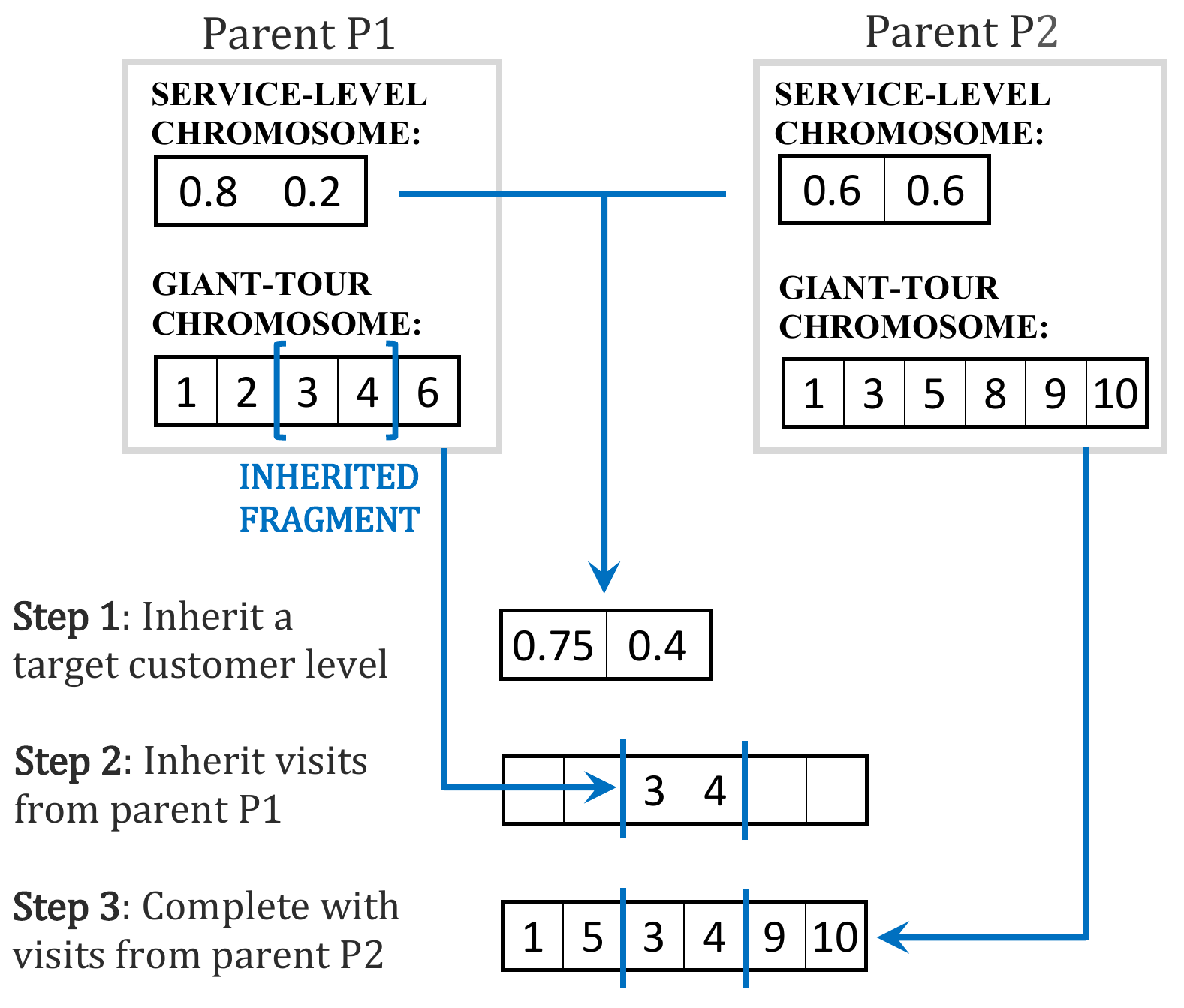}
 \caption{AOX crossover -- In this example, $\mathcal{V}_1 = \{1,2,3,4,5\}$, $\mathcal{V}_2 = \{6,7,8,9,10\}$ and $s_i = 1$ for all $i$.\label{fig:crossaox}}
\end{figure}

\paragraph{Selection and Crossover} To generate a new solution, the algorithm first selects two parents $P_1$ and~$P_2$ in the population via a binary tournament based on the \emph{biased fitness} measure described in Section \ref{section:pop}. 
A new offspring solution~$C$ is then obtained by crossover of $P_1$ and $P_2$.
For this purpose, we propose an {\em adapted order crossover} (AOX), which extends the well-known order crossover with the ability to transmit customer selection and visit sequence decision from both parents.
This crossover operator is illustrated in Figure \ref{fig:crossaox}.

In a first step, the crossover inherits the service level information from both parents. This is done by crossing the service level chromosomes of both parents using an \emph{extended intermediate recombination} \citep{Muhlenbein1993}, i.e., a target weight ratio $\alpha^{T}_k(C)$ is randomly chosen between $\alpha_k(P_1)$ and $\alpha_k(P_2)$ for each group~$k$, where $\alpha_k(P)$ is the weight ratio of group $k$ in individual $P$.

\myblue
{
Then, in a second step, the \myblue{giant-tour} chromosome of the child $C$ is initialized with the longest size among both parent and inherits, as in the \emph{order crossover} (OX), a fragment of $P_1$.

In the third and final step, the \myblue{giant-tour} chromosome of $C$ is completed by sweeping circularly the deliveries of $P_2$ and inheriting them, starting one index after the end of the fragment from~$P_1$. Each insertion of a visit $i$ of a group $k$ is done under the condition that $i$ does not already exist in $C$, and that the target service level $\alpha^{T}_k(C)$ has not yet been reached. To complete the representation, the service level chromosome of the child is finally derived from the giant-tour chromosome.}

\paragraph{Education}
The goal of the crossover operator was to generate new solutions which inherit common characteristics from both parents while introducing a significant level of randomness. As such, the crossover operator is not the main force which drives solution improvement, this role being assumed by a subsequent local search-based education procedure.

The local search (LS) is applied on the complete solution representation, including the visits to the depot. Therefore, the Split algorithm has to be run beforehand. The LS \myblue{uses} the same classical vehicle routing neighborhoods as in \cite{Vidal2012}: \textsc{2-Opt}, \textsc{2-Opt*}, \textsc{Swap}, \textsc{Relocate} as well as generalized \textsc{Swap} and \textsc{Relocate} involving two consecutive nodes, and limited to moves between close services. These classical neighborhoods only involve services which are present in the current solution. To also optimize the decision subset related to customer selections in the LS, we include three additional neighborhoods:
\begin{itemize}[nosep]
 \item[--] \textsc{Remove}: If $u$ is a visited customer, then remove $u$ from the solution.
 \item[--] \textsc{Add}: If $u$ is a visited customer and $v$ is a non-visited customer, add $v$ after~$u$.
 \item[--] \textsc{Replace:} If $u$ is a visited customer and $v$ is a non-visited customer, replace $u$ by $v$.
\end{itemize}
All neighborhoods are explored in a random order with a first improvement move acceptance policy. The LS stops when no improving move can be found in the entire neighborhood, and the resulting solution is converted back into a \myblue{giant-tour} and service level individual representation, which is inserted in the population.

\paragraph{Repair} Finally, it is possible for a solution to remain infeasible after education. When this situation happens, a \emph{Repair} operator is applied with $50\%$ chance. As in \citep{Vidal2012}, this operator simply consists in running the local search with higher penalty values, first by a factor of 10, and then 100, with the aim of converging towards a feasible solution.

\subsection{Population Management}
\label{section:pop}

The population is formed of two sub-populations, designed to host feasible and infeasible individuals, respectively. The algorithm starts by generating $4 \times \mu$ random initial solutions.
Each sub-population is then managed to contain at least $\mu$ and at most $\mu+\lambda$ individuals. Whenever a population reaches its maximum size, $\lambda$ individuals are eliminated to produce the next generation. This is done by iteratively selecting out either a \emph{clone} solution, with a distance of $0$ to another solution, or the worst solution according to a biased fitness function $\phi^\textsc{bias}(S)$ when there are no more clones.
The biased fitness measure evaluates every solution $S$ based on its cost $\phi^\textsc{cost}$ (Equation~\ref{solution-cost}) and its contribution to the sub-population diversity, defined as:
\begin{equation}
 \phi^\textsc{div}_P(S)=\frac{\sum_{S_2 \in N_P(S)}\delta(S, S_2)}{n^{\textsc{close}}},
\end{equation}
where $N_P(S)$ is the set of the $n^{\textsc{close}}$ individuals in the sub-population $P$ closest to $S$ with respect to a distance measure $\delta(S_1, S_2)$. \myblue{Our distance measure counts} the percentage of common edges between two solutions. Let $E(S_1)$ and $E(S_2)$ be the set of edges used in the solutions $S_1$ and $S_2$, then the distance is expressed as:
\begin{equation}
 \delta(S_1,S_2)=1-\frac{|E(S_1) \cap E(S_2)|}{|E(S_1) \cup E(S_2)|}.
\end{equation}

Finally, let $n^\textsc{elite}$ be a parameter controlling how elitist $\phi^\textsc{bias}(S)$ is, and let $R(S,P,f)$ be an application which returns the rank of an individual $S$ in the sub-population $P$ relatively to a measure~$f$. Then, the biased fitness of $S$ in $P$ is evaluated as in UHGS \citep{Vidal2012b}:
\begin{equation}\label{eq:biased_fitness}
 \phi^\textsc{bias}_P(S) = R(S,P,\phi^\textsc{cost}) + \Big(1 - \frac{n^\textsc{elite}}{|P|}\Big) R(S,P,\phi^\textsc{div}).
\end{equation}
\myblue{
Each distance measure can be evaluated in $\mathcal{O}(n)$ time. As such, when a new individual enters the population, its distance from every other solution can be computed in $\mathcal{O}(n |P|)$ time, and the computational complexity of updating the biased fitness measures is $\mathcal{O}(n|P| + |P| \log |P|)$.}

\paragraph{Adaptive Penalty Coefficients} The exploration of infeasible solutions contributes positively to the search if the ratio of feasible solutions is adequately controlled. To that extent, the penalty parameters are adapted to achieve a ratio of feasible solutions within a predefined interval $[\xi^\textsc{min},\xi^\textsc{max}]$.

We rely on $1+K$ penalty parameters for the VRP-SL: one for the capacity constraints, and~$K$ penalties for the service level constraints, one for each group. Let $\xi^c$ be the ratio of feasible individuals with respect to a constraint $c$, measured among the last 100 individuals generated by local search. In order to drive the search towards feasible solutions, every 100 iterations we update the penalty coefficient of constraint $c$ using the following rule:
\begin{equation}
w^c = 
\begin{cases}
w^c \times 1.2 & \hspace*{0.3cm} \text{ if } \xi^c \leq \xi^\textsc{min}, \\
 w^c \times 0.85 & \hspace*{0.3cm} \text{ if } \xi^c \geq \xi^\textsc{max},\\
 w^c & \hspace*{0.3cm} \text{ otherwise.}
\end{cases}
\end{equation}

Initially, all penalty parameters are set to $10$.
We aim to obtain around \myblue{25\%} feasible solutions after LS. Since all constraints need to be respected to obtain a globally feasible solution, we used $\xi^\textsc{min}=\myblue{0.15}^{\frac{1}{1+K}}$ and $\xi^\textsc{max}=\myblue{0.35}^{\frac{1}{1+K}}$ to achieve this goal.

\paragraph{Diversification procedure}
Finally, after each consecutive $It_\textsc{div}$ iterations without improvement of the best solution, we 
apply the same diversification procedure as in~\cite{Vidal2012} to only keep the best $\mu / 3$ individuals in each subpopulation and reintroduce new random initial solutions. This procedure complements the biased fitness function so as to avoid a premature convergence of the method due to the strong intensification of the LS.

\section{Experimental Analyses}
\label{section:experiments}

This section aims to
1) introduce a set of instances for the VRP-SL derived from classical vehicle routing instances,
\myblue{2) evaluate the performance of the proposed methods on the VRP-SL instances,
3) evaluate the performance of the metaheuristic and the branch-and-price on classical problems generalized by the VRP-SL, namely the VRPPFCC and CPTP, and finally}
4) examine the impact of some key parameters and design choices.
All algorithms were coded in C/C++.
\myblue{The exact algorithms were run on a single thread of a 3.07-GHz Intel Xeon CPU, and the metaheuristic was run on a single thread of a 3.4-GHz Intel Core i7 CPU. We used CPLEX 12.7 for the resolution of the compact formulation and for the linear programs in the B\&P algorithm.} \\

\noindent
\textbf{Benchmark instances for the VRP-SL.}
To compensate for the unavailability of benchmark instances for the VRP-SL, we derived two sets of instances from classical CVRP and prize-collecting VRP test sets:

\begin{itemize}[nosep,leftmargin=0.6cm]
\item[--] The first set (\textbf{S1}) has been generated from a subset of 26 instances from \cite{Augerat1995}, and includes between 31 and 80 vertices.
The profit of each client has been set to $h \times q_i$, where $h$ is a random variable uniformly generated in the interval $[0.75, 2.25]$.

\item[--] The second set (\textbf{S2}) has been derived from the 10 capacitated profitable tour instances of \cite{Archetti2013a}, with 51 to 200 vertices. The capacity of the vehicles has been set to 500, and the number of vehicles has been set to 
\begin{equation}
m = \left\lceil\frac{\sum_{k=1}^K (Q_k^\textsc{min} + Q_k^\textsc{max})}{2Q}\right\rceil,
\end{equation}
where $Q_k^\textsc{min}$ is the subset of services of group $k$ with smallest delivery quantity which allows to satisfy the service level requirements, and $Q_k^\textsc{max}$ is the sum of demands of all customers of this group. The original profits were multiplied by a factor of $0.5$ to obtain a good balance between profits and distance.
\end{itemize}

For each instance, we considered five configurations for the assignment of visits to groups: \mbox{\{1, 2R, 2C, 5R, 5C\}}. The first number corresponds to the number of groups, and the second letter, when applicable, corresponds to their distribution: random (R) or clustered (C). In all cases, the service levels for each group have been randomly selected in $\{0.45, 0.55, 0.65, 0.75, 0.85, 0.95,1\}$, and the weight of each customer, reflecting its importance for the service level constraint, coincides with the demand quantity (\mbox{$s_i$ = $q	_i$}). Overall, this leads to $26 \times 5 = 130$ instances of set~\textbf{S1}, and $10 \times 5 = 50$ instances of set~\textbf{S2}. For each instance, the distance values between customers are rounded to the nearest integer. For the sake of brevity, the results presented in the paper are aggregated per group of five instances. All instances and detailed results are available in the electronic companion of this paper, also available at \url{https://w1.cirrelt.ca/~vidalt/en/VRP-resources.html}.

\subsection{Exact solutions and lower bounds}
\label{results-exact}

\begin{table}[htbp]
\renewcommand{\arraystretch}{1.05}
\setlength{\tabcolsep}{0.35cm}
\vspace*{-0.7cm}
\centering
\rotatebox{90}{%
 \begin{varwidth}{1.1\textheight}
\scalebox{0.82}
{
\begin{tabular}{|r|rrr|rrrrrr|rrrrrr|}
\hline
&\multicolumn{3}{c|}{} & \multicolumn{6}{c|}{\textbf{Compact Formulation}} & \multicolumn{6}{c|}{\textbf{Branch-and-Price}} \\
\textbf{Data}&$\mathbf{n}$&$\mathbf{m}$&$\mathbf{n_\textsc{min}}$&\textbf{LB$_0$}&\textbf{T$_0$(s)}&\textbf{LB}&\textbf{Tree}&\textbf{Gap}&\textbf{T(s)}&\textbf{LB$_0$}&\textbf{T$_0$(s)}&\textbf{LB}&\textbf{Tree}&\textbf{Gap}&\textbf{T(s)}\\
\hline
A1 & 31 & 5 & 19.4 & 643.1 & 0.50 & 701.6 & 756k & 0.78 & 3423.93 & 685.9 & 0.13 & \textbf{707.4} & 319 & 0.00 & 33.82 \\
A2 & 32 & 6 & 17.6 & 610.0 & 0.37 & 652.7 & 650k & 0.77 & 2288.52 & 643.8 & 0.12 & \textbf{656.1} & 2k & 0.30 & 1443.06 \\
A3 & 35 & 5 & 21.8 & 627.1 & 0.40 & \textbf{666.2} & 407k & 0.00 & 1806.92 & 652.4 & 0.15 & \textbf{666.2} & 1k & 0.00 & 224.77 \\
A4 & 36 & 6 & 18.0 & 757.7 & 0.56 & 802.4 & 386k & 0.28 & 2482.49 & 792.6 & 0.13 & \textbf{805.0} & 103 & 0.00 & 5.93 \\
A5 & 38 & 5 & 18.8 & 629.3 & 0.55 & 667.7 & 584k & 0.25 & 3212.89 & 649.7 & 0.21 & \textbf{667.9} & 3k & 0.22 & 1861.53 \\
A6 & 43 & 7 & 25.4 & 739.9 & 0.88 & 772.7 & 414k & 1.94 & 4477.12 & 772.9 & 0.25 & \textbf{788.6} & 4k & 0.11 & 2890.31 \\
A7 & 44 & 7 & 28.2 & 956.9 & 1.66 & 989.3 & 560k & 3.25 & 7200.00 & 1006.1 & 0.31 & \textbf{1020.8} & 8k & 0.21 & 4118.13 \\
A8 & 47 & 7 & 26.6 & 874.5 & 1.40 & 905.9 & 670k & 4.47 & 7200.00 & 926.9 & 0.21 & \textbf{944.0} & 7k & 0.44 & 3394.61 \\
A9 & 53 & 7 & 34.2 & 993.4 & 1.94 & 1015.1 & 439k & 5.30 & 7200.00 & 1045.0 & 0.67 & \textbf{1066.0} & 8k & 0.64 & 4905.96 \\
A10 & 59 & 9 & 37.4 & 1136.6 & 3.18 & 1163.1 & 390k & 5.67 & 7200.00 & 1208.3 & 0.61 & \textbf{1226.3} & 14k & 0.60 & 5981.17 \\
A11 & 61 & 8 & 33.0 & 1033.1 & 3.62 & 1057.3 & 249k & 5.90 & 7200.00 & 1094.7 & 0.74 & \textbf{1112.1} & 7k & 1.04 & 4953.56 \\
A12 & 62 & 9 & 35.2 & 1055.1 & 2.54 & 1080.0 & 350k & 4.58 & 7200.00 & 1113.6 & 0.60 & \textbf{1128.3} & 18k & 0.35 & 6236.87 \\
A13 & 64 & 9 & 41.6 & 1030.2 & 2.17 & 1052.7 & 529k & 4.40 & 7200.00 & 1074.9 & 0.56 & \textbf{1090.4} & 11k & 1.00 & 5764.22 \\
A14 & 79 & 10 & 43.2 & 1432.9 & 6.00 & 1453.1 & 177k & 4.70 & 7200.00 & 1494.4 & 1.39 & \textbf{1508.6} & 7k & 1.09 & 7200.00 \\
B1 & 30 & 5 & 21.6 & 526.0 & 0.38 & \textbf{576.3} & 1786k & 1.42 & 5794.55 & 550.3 & 0.15 & 569.7 & 8k & 2.37 & 5762.20 \\
B2 & 34 & 5 & 16.8 & 695.4 & 0.26 & \textbf{736.0} & 1128k & 1.51 & 3476.90 & 709.6 & 0.31 & 733.9 & 6k & 1.91 & 5760.99 \\
B3 & 38 & 6 & 18.0 & 434.4 & 0.52 & 460.6 & 1673k & 7.24 & 7200.00 & 466.5 & 0.60 & \textbf{477.9} & 10k & 3.84 & 7200.00 \\
B4 & 42 & 6 & 22.0 & 618.2 & 0.83 & 650.7 & 766k & 2.29 & 6352.01 & 652.0 & 0.46 & \textbf{663.2} & 3k & 0.44 & 2705.18 \\
B5 & 44 & 5 & 22.8 & 548.8 & 0.82 & 587.0 & 1040k & 6.13 & 7200.00 & 606.2 & 0.77 & \textbf{614.8} & 7k & 1.68 & 7200.00 \\
B6 & 49 & 7 & 30.0 & 599.1 & 1.17 & 630.1 & 1049k & 5.74 & 7200.00 & 646.5 & 0.91 & \textbf{657.2} & 7k & 1.69 & 5140.14 \\
B7 & 50 & 7 & 35.4 & 863.7 & 0.58 & 908.5 & 1259k & 4.87 & 7200.00 & 916.5 & 0.91 & \textbf{933.3} & 6k & 2.23 & 5950.49 \\
B8 & 55 & 7 & 37.2 & 563.0 & 1.29 & 574.3 & 1117k & 9.76 & 7200.00 & 605.9 & 1.45 & \textbf{613.1} & 7k & 3.65 & 7200.00 \\
B9 & 56 & 9 & 35.2 & 1291.6 & 1.34 & 1316.1 & 724k & 4.32 & 7200.00 & 1332.6 & 0.63 & \textbf{1347.1} & 14k & 2.04 & 5761.09 \\
B10 & 63 & 9 & 34.0 & 726.2 & 1.97 & 761.4 & 731k & 6.22 & 7200.00 & 781.9 & 1.07 & \textbf{793.5} & 8k & 2.27 & 7200.00 \\
B11 & 66 & 10 & 42.0 & 878.6 & 1.85 & 907.8 & 787k & 4.54 & 7200.00 & 923.4 & 1.21 & \textbf{936.1} & 8k & 1.57 & 7200.00 \\
B12 & 77 & 10 & 43.0 & 1003.3 & 3.77 & 1020.3 & 450k & 6.49 & 7200.00 & 1060.6 & 1.91 & \textbf{1069.1} & 9k & 2.02 & 7200.00 \\
p03 & 100 & 3 & 53.8 & 542.1 & 3.48 & \textbf{554.4} & 176k & 0.11 & 2734.21 & 545.0 & 92.27 & 548.4 & 209 & 1.21 & 7200.00 \\
p06 & 50 & 2 & 30.0 & 370.4 & 0.44 & \textbf{387.0} & 1k & 0.00 & 5.26 & 374.1 & 11.61 & 381.0 & 418 & 1.54 & 7200.00 \\
p07 & 75 & 2.8 & 42.8 & 506.9 & 1.41 & \textbf{519.6} & 69k & 0.00 & 301.26 & 511.4 & 28.14 & 515.0 & 519 & 0.89 & 6489.66 \\
p08 & 100 & 3 & 55.2 & 531.3 & 3.33 & \textbf{544.3} & 109k & 0.09 & 1784.18 & 534.4 & 44.19 & 537.6 & 232 & 1.34 & 7200.00 \\
p09 & 150 & 4.6 & 88.2 & 669.4 & 17.11 & \textbf{684.5} & 2170k & 1.58 & 7200.00 & 681.9 & 306.57 & 684.0 & 79 & 1.64 & 7200.00 \\
p10 & 199 & 6 & 114.4 & 773.8 & 61.80 & 783.4 & 16k & 3.31 & 7200.00 & 798.2 & 699.11 & \textbf{799.6} & 45 & 1.32 & 7200.00 \\
p13 & 120 & 3 & 75.0 & 529.0 & 10.03 & 538.9 & 114k & 15.77 & 7200.00 & 587.3 & 837.90 & \textbf{588.9} & 27 & 7.96 & 7200.00 \\
p14 & 100 & 4 & 55.2 & 491.8 & 4.83 & 504.9 & 227k & 14.13 & 7200.00 & 560.3 & 90.53 & \textbf{568.5} & 575 & 3.31 & 7200.00 \\
p15 & 100 & 3.8 & 65.6 & 493.4 & 4.72 & 506.9 & 212k & 13.92 & 7200.00 & 565.4 & 108.80 & \textbf{572.1} & 557 & 2.85 & 7200.00 \\
p16 & 199 & 6.2 & 108.2 & 769.7 & 66.16 & 782.2 & 204k & 2.69 & 7200.00 & 792.7 & 571.70 & \textbf{794.0} & 37 & 1.22 & 7200.00 \\
\hline
\end{tabular}
}
 \caption{Results of the exact methods for the VRP-SL instances \label{tab:ex3a}}
 \end{varwidth}}
\end{table}

\myblue{
The two exact methods have been tested on each VRP-SL instance using a single core with a time limit of two hours. Since the VRP-SL instances were designed to be challenging for both exact and heuristic approaches, only a subset of the problems could be solved to optimality. To speed up the resolution, we used the best integer solution found by the metaheuristic as a warm start for both exact methods, hence limiting the size of the branch-and-bound tree.}

\myblue{
Table~\ref{tab:ex3a} presents average results for each group of instances. The first columns list the characteristics of each group, followed by the minimal number of visits $n_\textsc{min}$ required to satisfy service level constraints, the value LB$_0$ of the lower bound and the processing time T$_0$ at the root node, the value LB of the best overall lower bound, the number of nodes in the search tree, the final integrality gap and the total CPU time. For each instance, the best lower bound is highlighted in boldface. We do not indicate the best upper bound found by each method, since no improvement was found over the initial value obtained by the metaheuristic.}

 \begin{figure}[htbp]
 \centering
 \includegraphics[width=0.7\textwidth]{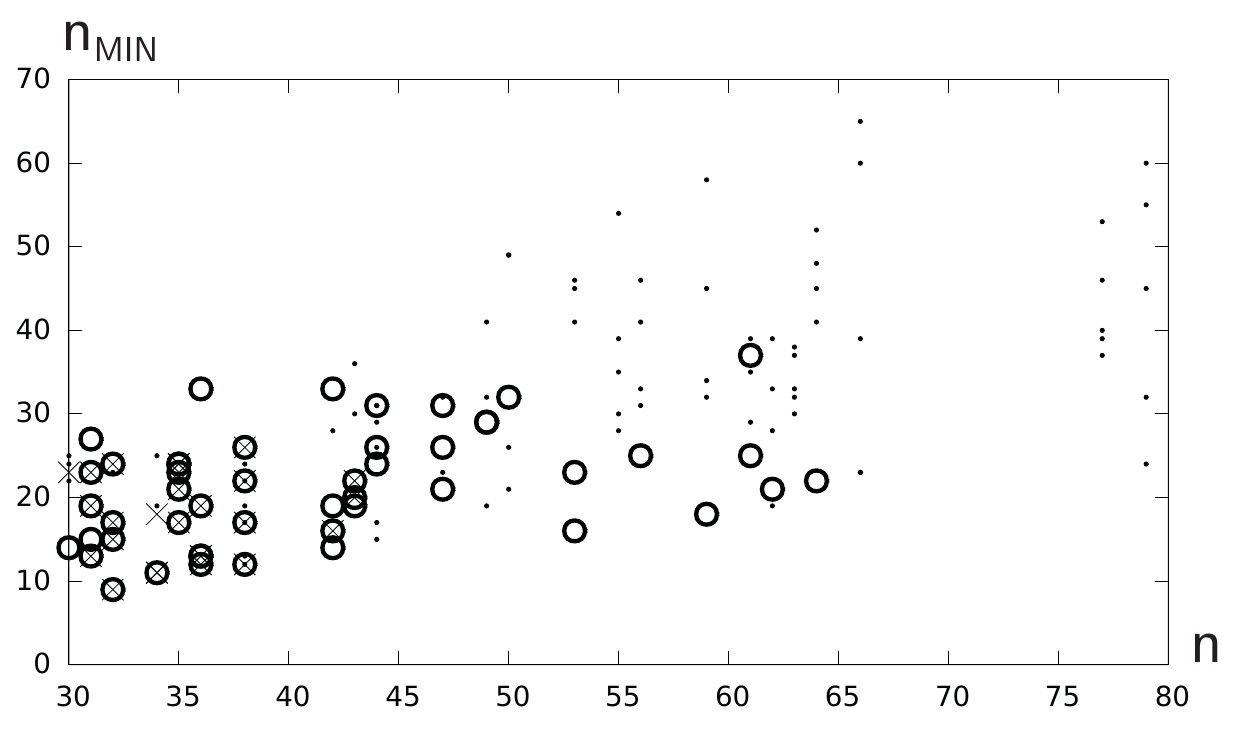}
 \caption{Instances of set \textbf{S1} currently solved to optimality \myblue{using the compact formulation ``$\times$'', via branch-and-price ``$\circ$'', or still open ``$\cdot$'',} as a function of the number of visits $n$ and the minimum number of deliveries $n_\textsc{min}$ \myblue{needed to satisfy the service level constraint}.}
 \label{fig:augerat_svrpp}
\end{figure}

\myblue{
As illustrated by the experiments, the branch-and-price algorithm outperformed the compact formulation-based method on most instances. More precisely, it performed better on all instances of the set \textbf{S1} with the exception of B1 and B2. For the set \textbf{S2}, the compact formulation performed better for the first five groups of instances and worse for the remaining ones. As expected, the branch-and-price algorithm is very effective when the solution includes short routes, while the compact formulation performs better on instances with few vehicles.
The total number of instances solved to optimality for configurations \{1, 2R, 2C, 5R, 5C\} were \{10, 9, 9, 9, 8\} for the compact formulation and \{15, 9, 9, 8, 9\} for the branch-and-price algorithm. These results show that the exact methods' performances are rather insensitive to the distribution of nodes in $V$, as well as the number of groups $K$. However, the complexity generally increases with $n$, $m$ and $n_\textsc{min}$, as visualized in Figure \ref{fig:augerat_svrpp}.
}

\subsection{Performance of the hybrid genetic algorithm}
\label{results-heuristic}

In this subsection, we report the results obtained with the proposed HGS on the VRP-SL instances. For each instance, the algorithm was run ten times with different seeds, using the same parameter setting and termination criterion as in \cite{Vidal2012} for the CVRP, that is $(n^\textsc{elite}, \mu, \lambda)=(10,25,40)$, $It_\textsc{max} = 2\times10^4$ and $It_\textsc{div} = 0.4 \times It_\textsc{max}$. Moreover, to accelerate the convergence, the parameters governing the population size have been halved when dealing with problem instances containing 200 or more services. Table \ref{tab:meta-VRPSL} reports, for each group of instances, \myblue{the worst, average and best solution quality of HGS over 10 runs  (Wor-10, Avg-10 and Best-10), the percentage gap between the average solutions and the best ones found (Gap$_\textsc{BKS}$), the percentage gap between the average solutions and the best lower bounds found by the exact methods (Gap$_\textsc{LB}$), the average CPU time in seconds, the best known solutions (BKS) and lower bounds (BKLB).}

\begin{table}[htbp]
\renewcommand{\arraystretch}{1.05}
\setlength{\tabcolsep}{0.25cm}
\centering
\scalebox{0.8}
{
\begin{tabular}{|r|cc|ccc@{\hspace*{0.4cm}}c@{\hspace*{0.4cm}}cc|cc|}
\hline
\textbf{Inst} &\textbf{n} & \textbf{m} & \textbf{Wor-10} &  \textbf{Avg-10} & \textbf{Best-10} &  \  \textbf{Gap$_\textsc{BKS}$} &  \textbf{Gap$_\textsc{LB}$}  & \textbf{T(s)}  & \textbf{BKS} & \textbf{BKLB} \\
\hline
A1  & 31 & 5 & 707.4 & 707.40 & 707.4 & 0.00 & 0.00 & 9.06 & 707.4 & 707.4 \\ 
A2  & 32 & 6 & 658.2 & 658.20 & 658.2 & 0.00 & 0.32 & 8.72 & 658.2 & 656.1 \\ 
A3  & 35 & 5 & 666.2 & 666.20 & 666.2 & 0.00 & 0.00 & 9.70 & 666.2 & 666.2 \\ 
A4  & 36 & 6 & 805.0 & 805.00 & 805.0 & 0.00 & 0.00 & 9.51 & 805.0 & 805.0 \\ 
A5  & 38 & 5 & 669.4 & 669.40 & 669.4 & 0.00 & 0.23 & 10.18 & 669.4 & 667.9 \\ 
A6  & 43 & 6 & 789.6 & 789.60 & 789.6 & 0.00 & 0.13 & 12.78 & 789.6 & 788.6 \\ 
A7  & 44 & 7 & 1023.0 & 1023.00 & 1023.0 & 0.00 & 0.22 & 13.90 & 1023.0 & 1020.8 \\ 
A8  & 47 & 7 & 948.2 & 948.20 & 948.2 & 0.00 & 0.44 & 12.99 & 948.2 & 944.0 \\ 
A9  & 53 & 7 & 1073.2 & 1073.20 & 1073.2 & 0.00 & 0.68 & 18.49 & 1073.2 & 1066.0 \\ 
A10  & 59 & 9 & 1234.8 & 1234.50 & 1234.4 & 0.01 & 0.67 & 21.21 & 1234.4 & 1226.3 \\ 
A11  & 61 & 8 & 1124.8 & 1124.10 & 1124.0 & 0.01 & 1.08 & 20.69 & 1124.0 & 1112.1 \\ 
A12  & 62 & 10 & 1132.6 & 1132.44 & 1132.4 & 0.00 & 0.36 & 19.54 & 1132.4 & 1128.3 \\ 
A13  & 64 & 9 & 1101.8 & 1101.70 & 1101.6 & 0.01 & 1.04 & 22.64 & 1101.6 & 1090.4 \\ 
A14  & 79 & 10 & 1526.2 & 1525.30 & 1525.2 & 0.01 & 1.11 & 32.05 & 1525.2 & 1508.6 \\ 
B1  & 30 & 5 & 584.8 & 584.80 & 584.8 & 0.00 & 0.57 & 9.53 & 584.8 & 581.5 \\ 
B2  & 34 & 5 & 748.2 & 748.20 & 748.2 & 0.00 & 1.10 & 8.34 & 748.2 & 740.0 \\ 
B3  & 38 & 5 & 497.2 & 497.20 & 497.2 & 0.00 & 3.98 & 10.24 & 497.2 & 478.2 \\ 
B4  & 42 & 6 & 666.2 & 666.20 & 666.2 & 0.00 & 0.46 & 12.00 & 666.2 & 663.2 \\ 
B5  & 44 & 5 & 625.0 & 625.00 & 625.0 & 0.00 & 1.66 & 11.71 & 625.0 & 614.8 \\ 
B6  & 49 & 7 & 668.8 & 668.80 & 668.8 & 0.00 & 1.77 & 14.25 & 668.8 & 657.2 \\ 
B7  & 50 & 7 & 956.2 & 956.20 & 956.2 & 0.00 & 2.38 & 15.13 & 956.2 & 934.0 \\ 
B8  & 55 & 7 & 637.4 & 637.40 & 637.4 & 0.00 & 3.97 & 18.26 & 637.4 & 613.1 \\ 
B9  & 56 & 9 & 1378.2 & 1376.72 & 1376.2 & 0.04 & 2.03 & 22.17 & 1376.2 & 1349.4 \\ 
B10  & 63 & 9 & 812.0 & 812.00 & 812.0 & 0.00 & 2.33 & 20.20 & 812.0 & 793.5 \\ 
B11  & 66 & 10 & 951.8 & 951.56 & 951.4 & 0.02 & 1.65 & 24.15 & 951.4 & 936.1 \\ 
B12  & 77 & 10 & 1092.6 & 1091.94 & 1091.4 & 0.05 & 2.13 & 28.92 & 1091.4 & 1069.1 \\ 
p03  & 100 & -- & 555.0 & 555.00 & 555.0 & 0.00 & 0.10 & 22.67 & 555.0 & 554.4 \\ 
p06  & 50 & -- & 387.0 & 387.00 & 387.0 & 0.00 & 0.00 & 10.74 & 387.0 & 387.0 \\ 
p07  & 75 & -- & 519.6 & 519.60 & 519.6 & 0.00 & 0.00 & 17.51 & 519.6 & 519.6 \\ 
p08  & 100 & -- & 544.8 & 544.80 & 544.8 & 0.00 & 0.09 & 22.50 & 544.8 & 544.3 \\ 
p09  & 150 & -- & 695.6 & 695.60 & 695.6 & 0.00 & 1.50 & 46.12 & 695.6 & 685.3 \\ 
p10  & 199 & -- & 810.8 & 810.44 & 810.2 & 0.03 & 1.36 & 86.30 & 810.2 & 799.6 \\ 
p13  & 120 & -- & 640.8 & 640.20 & 639.8 & 0.06 & 8.70 & 44.05 & 639.8 & 588.9 \\ 
p14  & 100 & -- & 588.0 & 588.00 & 588.0 & 0.00 & 3.43 & 27.04 & 588.0 & 568.5 \\ 
p15  & 100 & -- & 588.8 & 588.80 & 588.8 & 0.00 & 2.93 & 26.78 & 588.8 & 572.1 \\ 
p16  & 199 & -- & 805.0 & 804.24 & 803.8 & 0.05 & 1.29 & 92.85 & 803.8 & 794.0 \\ 
\hline
\multicolumn{3}{|r|}{All}  &&& & 0.01 & 1.38 & 22.58 & & \\
\hline
\end{tabular}
}
\caption{Performance of the HGS on the VRP-SL instance sets \label{tab:meta-VRPSL}}
\end{table}

In the absence of results from previously published heuristics, we consider \myblue{three main} indicators of method performance: its ability to reach known optimal solutions, the stability of the solution quality over several runs, illustrated by the gap between the average solution quality and the BKS, and \myblue{the deviation from the lower bounds produced by the mathematical programming algorithms}. As observed in our experiments, 
HGS finds all of the \myblue{70}
known optimal solutions on all test runs. 
The method also returns solutions of consistent high quality: it found for \myblue{25}/36 groups of instances the value of the best known solution on all runs of all instances. \myblue{The percentage gaps between the average and best known solutions are close to zero ($0.01\%$ overall), and the instances with fewer groups appear to be generally easier to solve. The gaps to the lower bounds are also small ($1.38\%$ in average), and thus the solutions of HGS are guaranteed to be close to the optima. This gap is more likely to be due to the quality of the lower bounds, since over all 2h-runs of the exact approaches not a single best solution of the metaheuristic was improved.}
Finally, the average CPU time never exceeds \myblue{93 seconds}, the worst case being observed on problem \myblue{p16} with 199 services.

\subsection{Comparative analyses on key subproblems}

As discussed in Section \ref{section:Litt}, the VRP-SL generalizes several important problem classes.
Two such problems in particular, the VRPPFCC and CPTP, have been the focus of a wide literature, opening the way to some comparative performance analyses.

\paragraph{Experiments on the VRPPFCC} We rely on the two sets of instances from \cite{Bolduc2008}. Set CE includes up to 199 customers, while Set G includes larger instances with up to 483 customers. \myblue{For these instances, the convention is to compute all distances with double precision and report the final solution with two digits. Our heuristic method is compared to the two best current metaheuristics in the literature: the UGHS proposed by \cite{Vidal2014}, which uses an \emph{exhaustive} solution representation (with all visits) with a route-evaluation operator in charge of customer selections, and the recent AVNS of \cite{Huijink2014}. To the best of our knowledge there are no known exact results for this problem in the literature.} The results of our methods are displayed in Table \ref{tab:ex7}. The column BKS reports the best known solutions found in the literature before this paper. Avg-10 and Best-10 represent the average and best solutions found by the \myblue{metaheuristic} over ten runs. For each instance, the best \myblue{heuristic} method in highlighted in boldface. \myblue{For the branch-and-price results, columns LB$_0$ and LB follow the same convention as Table~\ref{tab:ex3a}.}

\begin{table}[htbp]
\renewcommand{\arraystretch}{1.1}
\setlength{\tabcolsep}{0.2cm}
\centering
\vspace*{-0.95cm}
\hspace*{-0.5cm}
\rotatebox{90}{%
 \begin{varwidth}{1.1\textheight}
\scalebox{0.8}
{
\begin{tabular}{|rrr|rrr|rrr|rrrH|rrr|r|}
\hline
 & & &\multicolumn{3}{c|}{\textbf{UGHS -- \cite{Vidal2014}}} &  \multicolumn{3}{c|}{\textbf{AVNS -- \cite{Huijink2014}}} & \multicolumn{4}{c|}{\textbf{This paper -- HGS}} & \multicolumn{3}{c|}{\textbf{This paper -- B\&P}} & \\
Instance & $n$ & $m$ & Avg-10 & Best-10 & T(s)  & Avg-10 & Best-10 & T(s) & Avg-10 & Best-10 & T(s) & T$^*$(s) & LB$_0$ & LB & T(s) & BKS\\
\hline
CE-01 & 50 & 4 & 1119.66 & \textbf{1119.47} & 38.5 & \textbf{1119.47} & \textbf{1119.47} & 47.4 & \textbf{1119.47} & \textbf{1119.47} & 13.01 & 0.68 & 1112.17 & 1119.47 & 154.30 & 1119.47 \\
CE-02 & 75 & 9 & 1815.63 & \textbf{1814.52} & 59.6 & 1817.06 & \textbf{1814.52} & 95.8 & \textbf{1814.52} & \textbf{1814.52} & 24.00 & 4.49 & 1795.13 & 1807.79 & 7200.00 & 1814.52 \\
CE-03 & 100 & 6 & \textbf{1922.88} & \textbf{1919.05} & 476.8 & 1925.08 & \textbf{1919.05} & 236.7 & 1923.92 & \textbf{1919.05} & 45.16 & 20.49 & 1896.81 & 1904.30 & 7200.00 & 1919.05 \\
CE-04 & 150 & 9 & \textbf{2509.82} & \textbf{2505.39} & 934.7 & 2518.14 & 2509.81 & 642.7 & 2510.35 & \textbf{2505.39} & 98.95 & 61.96 & 2478.50 & 2483.50 & 7200.00 & 2505.39 \\
CE-05 & 199 & 13 & 3095.58 & \textbf{3081.59} & 1289.3 & 3101.40 & 3090.49 & 1437.2 & \textbf{3094.78} & 3084.40 & 213.69 & 159.57 & 3051.13 & 3053.92 & 7200.00 & 3081.59 \\
CE-06 & 50 & 4 & \textbf{1207.47} & \textbf{1207.47} & 38.5 & \textbf{1207.47} & \textbf{1207.47} & 44.2 & \textbf{1207.47} & \textbf{1207.47} & 12.42 & 0.59 & 1199.81 & 1207.47 & 196.41 & 1207.47 \\
CE-07 & 75 & 9 & 2012.33 & 2006.52 & 73.2 & 2009.58 & \textbf{2004.53} & 98.6 & \textbf{2006.32} & \textbf{2004.53} & 25.71 & 5.65 & 1984.80 & 1997.90 & 7200.00 & 2004.53 \\
CE-08 & 100 & 6 & \textbf{2057.57} & \textbf{2052.05} & 500.3 & 2059.56 & \textbf{2052.05} & 239.6 & 2058.18 & \textbf{2052.05} & 48.30 & 24.22 & 2029.89 & 2037.66 & 7200.00 & 2052.05 \\
CE-09 & 150 & 10 & 2428.19 & 2425.32 & 1241.0 & 2426.35 & 2420.71 & 953.9 & \textbf{2422.88} & \textbf{2420.57} & 123.43 & 82.62 & 2389.86 & 2394.60 & 7200.00 & 2419.84 \\
CE-10 & 199 & 13 & 3387.12 & 3381.67 & 1229.9 & 3388.22 & \textbf{3373.84} & 1704.1 & \textbf{3384.03} & 3374.02 & 238.95 & 185.7 & 3340.29 & 3343.29 & 7200.00 & 3373.84 \\
CE-11 & 120 & 6 & 2331.13 & \textbf{2330.94} & 1202.9 & \textbf{2330.94} & \textbf{2330.94} & 399.0 & 2330.95 & \textbf{2330.94} & 71.94 & 43.69 & 2317.21 & 2321.26 & 7200.00 & 2330.94 \\
CE-12 & 100 & 8 & 1953.13 & \textbf{1952.86} & 150.8 & \textbf{1952.86} & \textbf{1952.86} & 183.6 & 1953.27 & \textbf{1952.86} & 34.28 & 8.66 & 1939.32 & 1946.63 & 7200.00 & 1952.86 \\
CE-13 & 120 & 6 & 2859.07 & \textbf{2858.83} & 1184.9 & \textbf{2858.92} & \textbf{2858.83} & 437.8 & 2858.94 & \textbf{2858.83} & 65.11 & 35.25 & 2843.90 & 2847.50 & 7200.00 & 2858.83 \\
CE-14 & 100 & 7 & \textbf{2213.02} & \textbf{2213.02} & 189.8 & 2213.78 & \textbf{2213.02} & 185.2 & 2213.41 & \textbf{2213.02} & 33.11 & 9.94 & 2193.70 & 2202.21 & 7200.00 & 2213.02 \\
G-01 & 240 & 7 & 14151.51 & 14131.18 & 2405.9 & 14163.43 & 14129.48 & 2433.9 & \textbf{14135.93} & \textbf{14112.39} & 210.40 & 147.8 & 14024.30 & 14038.00 & 7200.00 & 14111.95 \\
G-02 & 320 & 8 & 19190.77 & 19166.58 & 2409.9 & 19254.23 & 19140.69 & 6630.0 & \textbf{19134.06} & \underline{\textbf{19125.35}} & 302.27 & 215.28 & 18974.40 & 18978.70 & 7200.00 & 19140.69 \\
G-03 & 400 & 8 & 24588.29 & 24409.02 & 2418.1 & 24566.00 & 24406.67 & 11062.8 & \textbf{24372.81} & \underline{\textbf{24323.64}} & 431.70 & 320.8 & --- & --- & 7200.00 & 24368.29 \\
G-04 & 480 & 8 & 34517.47 & 34362.8 & 2421.1 & 34425.00 & 34231.56 & 15875.4 & \textbf{34116.56} & \underline{\textbf{34076.38}} & 401.61 & 287.92 & --- & --- & 7200.00 & 34231.56 \\
G-05 & 200 & 4 & 14296.07 & \textbf{14223.63} & 2408.0 & \textbf{14261.06} & 14229.50 & 1591.3 & 14268.87 & \textbf{14223.63} & 71.98 & 25.14 & --- & --- & 7200.00 & 14223.63 \\
G-06 & 280 & 5 & 21488.29 & 21396.60 & 2411.4 & 21440.38 & \textbf{21357.16} & 4337.2 & \textbf{21376.19} & \textbf{21357.16} & 125.28 & 69.81 & --- & --- & 7200.00 & 21357.16 \\
G-07 & 360 & 7 & 23463.05 & 23373.38 & 2414.9 & 23440.51 & 23263.22 & 7585.3 & \textbf{23285.33} & \underline{\textbf{23233.97}} & 354.48 & 256.37 & --- & --- & 7200.00 & 23263.22 \\
G-08 & 440 & 8 & 29918.06 & 29823.18 & 2415.6 & 29864.19 & 29657.38 & 12316.5 & \textbf{29603.54} & \underline{\textbf{29556.13}} & 415.88 & 311.13 & --- & --- & 7200.00 & 29657.38 \\
G-09 & 255 & 11 & 1332.63 & 1328.65 & 2323.4 & 1325.79 & \textbf{1320.29} & 3852 & \textbf{1324.06} & 1320.86 & 335.90 & 256.94 & 1308.14 & 1308.59 & 7200.00 & 1319.72 \\
G-10 & 323 & 13 & 1603.82 & 1597.61 & 2342.3 & 1592.14 & 1588.05 & 6922.9 & \textbf{1587.42} & \textbf{1584.73} & 647.96 & 538.61 & 1567.66 & 1568.01 & 7200.00 & 1583.50 \\
G-11 & 399 & 14 & 2192.68 & 2182.01 & 2405.2 & 2172.45 & 2163.50 & 12303.3 & \textbf{2166.52} & \textbf{2160.43} & 861.43 & 749.98 & 2137.91 & 2138.06 & 7200.00 & 2159.78 \\
G-12 & 483 & 15 & 2529.84 & 2522.64 & 2407.2 & 2493.94 & 2483.06 & 19555.0 & \textbf{2486.56} & \textbf{2481.76} & 956.42 & 825.02 & 2455.51 & 2455.64 & 7200.00 & 2479.62 \\
G-13 & 252 & 21 & \textbf{2261.50} & \textbf{2258.02} & 1412.7 & 2264.92 & 2261.66 & 2260.2 & 2266.41 & 2263.48 & 325.29 & 265.07 & 2232.70 & 2232.87 & 7200.00 & 2258.02 \\
G-14 & 320 & 23 & \textbf{2687.50} & \textbf{2683.73} & 1935.7 & 2689.99 & 2684.66 & 4838.5 & 2696.11 & 2689.67 & 682.99 & 596.26 & 2654.94 & 2655.23 & 7200.00 & 2683.73 \\
G-15 & 396 & 26 & \textbf{3152.00} & 3145.11 & 2301.4 & 3156.84 & 3150.67 & 8198.0 & 3153.76 & \underline{\textbf{3144.51}} & 858.02 & 761.46 & 3107.12 & 3107.42 & 7200.00 & 3145.11 \\
G-16 & 480 & 29 & \textbf{3632.04} & \textbf{3620.71} & 2450.1 & 3633.76 & 3624.56 & 12741.7 & 3635.54 & 3624.91 & 954.65 & 828.38 & 3578.85 & 3579.04 & 7200.00 & 3620.71 \\
G-17 & 240 & 18 & 1671.72 & \textbf{1666.31} & 1805.3 & 1672.81 & \textbf{1666.31} & 1878 & \textbf{1667.81} & \textbf{1666.31} & 160.84 & 104.76 & 1666.31 & 1666.31 & 19.69 & 1666.31 \\
G-18 & 300 & 22 & \textbf{2733.12} & \textbf{2730.55} & 2035.0 & 2734.86 & 2731.28 & 3916.6 & 2733.49 & 2730.97 & 459.92 & 385.05 & 2717.24 & 2717.95 & 7200.00 & 2730.55 \\
G-19 & 360 & 26 & 3504.26 & 3497.20 & 1989.5 & 3499.55 & 3494.28 & 5225.4 & \textbf{3495.68} & \underline{\textbf{3491.24}} & 839.89 & 749.46 & 3477.87 & 3478.79 & 7200.00 & 3494.27 \\
G-20 & 420 & 31 & 4319.37 & 4312.45 & 2523 & 4314.48 & 4307.63 & 8081.6 & \textbf{4307.61} & \underline{\textbf{4300.77}} & 1115.09 & 1059.51 & 4277.33 & 4278.20 & 7200.00 & 4306.85 \\
\hline
\multicolumn{3}{|c|}{Avg. Gap(\%)}   & 0.445 & 0.210 &  & 0.345 & 0.058 &  & 0.141 & -0.012 &  &  & --- & --- & & \\ 
\multicolumn{3}{|c|}{Avg. T(s)}   &  &  & 1583.70 &  &  & 4656.22 &  &  &  & 276.41 & & & 6575.60 & \\ 
\multicolumn{3}{|c|}{CPU}   & \multicolumn{3}{c|}{Xe 3.07GHz}  & \multicolumn{3}{c|}{Intel i5 2.6GHz} & \multicolumn{4}{c|}{Intel i7 3.4GHz} & \multicolumn{3}{c|}{Intel i7 3.4GHz} &    \\ 
 \hline
\end{tabular}
}
\caption{Performance of the HGS on VRPPFCC benchmark instances, in comparison with the current state-of-the-art algorithms \label{tab:ex7}}
 \end{varwidth}}
\end{table}

In these experiments, the proposed \myblue{metaheuristic} appears to outperform previous methods in terms of solution quality, with an average gap of \myblue{0.141\%}, in comparison to 0.445\% for UHGS with the exhaustive solution representation, and 0.345\% for AVNS.  During these tests, eight new best known solutions (BKS) have been found, as underlined in the table. Finally, the average CPU time is \myblue{markedly} faster than previous approaches which were run on processors of a similar generation, with \myblue{340} seconds on average, compared to 1584 and 4656 seconds for the other methods. \myblue{The proposed branch-and-price algorithm is the first in the literature to report optimal VRPPFCC solutions for three instances: CE-01, CE-06 and G-17. On the other hand, on six instances, the algorithm was not able to solve the root node within a time limit of two hours. This is due to their size (up to 480 customers), which can be considered very large for the current exact methods. Finally, the average integrality gap calculated over the tractable instances was 0.667\%, confirming the good performance of the approach.}

\paragraph{Experiments on the CPTP} For this problem, we compare the proposed HGS with the previous two best methods in the literature, the UHGS with \emph{exhaustive} solution representation and the multi-start ILS of \cite{Vidal2014}. \myblue{We also compare the proposed branch-and-price algorithm with the exact approach of \cite{Archetti2013a}, which generated, to this date, the best bounds for the CPTP. In that work, the authors proposed a branch-and-price algorithm using a \textit{q}-route relaxation, and reported detailed results with and without a primal heuristic.}
\myblue{We rely on the ten test cases of~\cite{Archetti2008d}, each case being used to produce $12$ instances with} a different fleet size and vehicle capacity, for a total of $120$ instances. \myblue{Tables~\ref{tab:summaryCPTP} and \ref{tab:summaryCPTP-BP} present a summary of our computational experiments on these instances,} using the same conventions as previously. Each line in the table corresponds to an average measure over a group of $12$ instances.

\begin{table}[htbp]
\renewcommand{\arraystretch}{1.1}
\hspace*{-0.4cm}
\centering
\setlength{\tabcolsep}{0.4cm}
\scalebox{0.84}
{
\setlength{\tabcolsep}{5pt}
\begin{tabular}{|rr|rrr|rrr|rrrH|r|}
\hline
&  &\multicolumn{3}{c|}{\textbf{UGHS -- \cite{Vidal2014}}} &  \multicolumn{3}{c|}{\textbf{ILS -- \cite{Vidal2014}}} & \multicolumn{4}{c|}{\textbf{This paper -- HGS}} & \\
Instance & $n$ & Avg-10 & Best-10 & T(s)  & Avg-10 & Best-10 & T(s) & Avg-10 & Best-10 & T(s) & T$^*$(s) & BKS \\
\hline
p03 & 100 & {\bf 254.07} & {\bf 254.07} & 215.82 & 253.99 & {\bf 254.07} & 191.63 & {\bf 254.07} & {\bf 254.07} & 12.32 & 2.43 & 254.07 \\
p06 & 50 & {\bf 129.13} & {\bf 129.13} & 30.78 & 129.09 & 129.11 & 23.35 & {\bf 129.13} &{\bf 129.13} & 7.18 & 0.16 & 129.13 \\
p07 & 75 & 192.51 & {\bf 192.56} & 103.94 & 192.37 & {\bf 192.56} & 85.80 & {\bf 192.52} & {\bf 192.56} & 8.36 & 0.57 & 192.56 \\
p08 & 100 & {\bf 254.07} & {\bf 254.07} & 216.28 & 253.93 & {\bf 254.07} & 191.49 & 253.90 & {\bf 254.07} & 12.17 & 2.15 & 254.07 \\
p09 & 150 & 319.61 & {\bf 319.72} & 295.34 & 319.23 & {\bf 319.72} & 289.75 & {\bf 319.67} & {\bf 319.72} & 18.44 & 5.42 & 319.72 \\
p10 & 199 & 387.87 & 388.41 & 303.55 & 386.67 & 387.73 & 309.23 & {\bf 388.79} & \underline{{\bf 388.85}} & 19.90 & 4.56 & 388.41 \\
p13 & 120 & {\bf 180.20} & {\bf 180.39} & 235.98 & 178.33 & 180.32 & 255.97 & 180.08 & 180.32 & 15.15 & 1.74 & 180.39 \\
p14 & 100 & {\bf 246.24} & {\bf 246.24} & 116.26 & 246.23 & {\bf 246.24} & 111.83 &{\bf 246.24} & {\bf 246.24} & 10.78 & 0.76 & 246.24 \\
p15 & 150 & 327.99 & 328.36 & 295.02 & 326.91 & 327.81 & 289.85 & {\bf 328.28} & \underline{{\bf 328.37}} & 15.24 & 2.10 & 328.36 \\
p16 & 199 & 393.09 & 393.75 & 303.88 & 392.21 & 393.32 & 309.74 & {\bf 394.03} & \underline{{\bf 394.04}} & 21.16 & 5.41 & 393.75 \\
\hline
\multicolumn{2}{|c|}{Avg. Gap(\%)}& 0.029 & 0.000 && 0.172 & 0.014 && 0.012 & -0.002 &&&\\
\multicolumn{2}{|c|}{Avg. T(s)}&&& 211.68 &&& 205.86 &&& 14.07 &&\\
\multicolumn{2}{|c|}{CPU} & \multicolumn{3}{c|}{Xe 3.07GHz}  & \multicolumn{3}{c|}{Xe 3.07GHz} & \multicolumn{4}{c|}{Intel i7 3.4GHz} & \\
\hline
\end{tabular}
}
\caption{Performance of HGS on the CPTP benchmark instances \label{tab:summaryCPTP}}
\end{table}

\begin{table}[htbp]
\renewcommand{\arraystretch}{1.1}
\hspace*{-0.4cm}
\centering
\setlength{\tabcolsep}{0.4cm}
\scalebox{0.84}
{
\setlength{\tabcolsep}{5pt}
\begin{tabular}{|rr|rrr|rrr|rrrr|r|}
\hline
&  &\multicolumn{3}{c|}{\textbf{B\&P1 -- \cite{Archetti2013a}}} &  \multicolumn{3}{c|}{\textbf{B\&P2 -- \cite{Archetti2013a}}} & \multicolumn{4}{c|}{\textbf{This paper -- B\&P}} & \\
Instance & $n$ & UB & LB & T(s) & UB & LB & T(s) & UB$_0$ & UB & LB & T(s) & BKUB \\
\hline
p03 & 100 & 256.90 &  {\bf 254.07} & 1135.31 & 256.82 & 147.01 & 1135.38 & 258.17 & {\bf \underline{255.36}} &   {\bf 254.07} & 1116.00 & 256.82 \\
p06 & 50 & 129.75 &  {\bf 129.13} & 625.08 & 129.64 & 118.09 & 602.38 & 131.77 & {\bf \underline{129.39}} &  {\bf 129.13} & 287.89 & 129.64 \\
p07 & 75 & 193.18 &  {\bf 192.56} & 622.38 & 193.11 & 163.49 & 613.92 & 194.21 & {\bf \underline{192.82}} &  {\bf 192.56} & 285.06 & 193.11 \\
p08 & 100 & 256.89 &  {\bf 254.07} & 1134.31 & 256.82 & 147.01 & 1135.00 & 258.17 & {\bf \underline{255.36}} & {\bf 254.07} & 1116.05 & 256.82 \\
p09 & 150 & 324.25 & 316.33 & 1163.69 & 324.20 & 160.90 & 1163.85 & 323.66 & {\bf \underline{320.83}} &  {\bf 319.72} & 1123.65 & 324.20 \\
p10 & 199 & 392.03 & 377.01 & 1117.77 & 391.66 & 176.43 & 1117.77 & 392.07 & {\bf \underline{390.43}} &  {\bf 388.85} & 596.37 & 391.66 \\
p13 & 120 & 191.74 & 167.78 & 1916.08 & 191.70 & 116.06 & 1916.77 & 186.72 & {\bf \underline{183.93}} & {\bf 180.32} & 1216.18 & 191.70 \\
p14 & 100 & 255.47 & 237.02 & 1110.31 & 255.41 & 116.24 & 1110.31 & 248.67 & {\bf \underline{247.20}} &  {\bf 246.24} & 556.22 & 255.40 \\
p15 & 150 & 331.98 & 324.19 & 893.54 & 331.93 & 190.18 & 893.77 & 332.36 & {\bf \underline{330.67}} &  {\bf 328.37} & 603.52 & 331.93 \\
p16 & 199 & 398.57 & 379.03 & 1118.46 & 398.27 & 178.90 & 1118.38 & 398.58 & {\bf \underline{396.47}} &  {\bf 394.04} & 941.83 & 398.27 \\
\hline
\multicolumn{2}{|c|}{Average} & 273.07 & 263.12 & 1083.69 & 272.96 & 151.43 & 1080.75 & 272.44 & 270.25 & 268.74 & 784.28 & 272.96 \\
\multicolumn{2}{|c|}{CPU} & \multicolumn{3}{c|}{Xe 2.26GHz}  & \multicolumn{3}{c|}{Xe 2.26GHz} & \multicolumn{4}{c|}{Intel i7 3.4GHz} & \\
\hline
\end{tabular}
}
\caption{Performance of the proposed branch-and-price algorithm on the CPTP benchmark instances \label{tab:summaryCPTP-BP}}
\end{table}

These instances are generally smaller, with $50$ to $199$ service locations, and thus the performance differences between state-of-the-art \myblue{heuristics} are less marked. Still, we observe that the proposed HGS obtains average solutions of similar or better quality (\myblue{0.012\%} gap compared to 0.029\% and 0.172\% gap) in a fraction of the CPU time of previous algorithms (\myblue{$14$}~seconds in average, compared to $211$ or $205$ seconds). \myblue{Three previous BKS were improved, leading to an average gap of $-0.002\%$ for the best solution quality of 10 runs.} 
\myblue{For the branch-and-price algorithm, we set the time limit to 3600 seconds in order to make a fair comparison with~\cite{Archetti2013a}. The proposed B\&P found similar or better solutions for all instances tested. It improved the bounds for 37 instances, with an average improvement of 0.629\%, and proved optimality for 103 instances, including ten new optimality certificates.}

%

\subsection{Sensitivity analyses}

Finally, this section reports additional sensitivity analyses on the impact of key components of the proposed HGS. For this purpose, we compare the results of the standard method described in Section \ref{section:HGA} against several alternative configurations obtained by modifying one operator, design choice, or group of parameters:\vspace*{0.2cm}

\begin{itemize}[nosep,leftmargin=2cm]
\item[\textbf{EOX} --] \myblue{An edge-recombination crossover (EOX) is used. As described in \citep{Whitley1989}, this crossover maintains, for each vertex $i$, an \emph{adjacency list} of non-visited vertices adjacent to $i$ in at least one parent. After a random choice for the first vertex, EOX iteratively inserts the adjacent vertex with the shortest adjacency list. Ties are broken randomly, and a random vertex is chosen whenever the adjacency list is empty. As in our adaptation of OX, target service levels are inherited for the groups from the parents, and any vertex belonging to a group for which the target service level has already been attained is eliminated from  the adjacency lists.}
\item[\textbf{No SL} --] Service levels are not used to filter service insertions in the crossover.
\item[\textbf{No INF} --] All penalty coefficients are set to a large value to avoid infeasibility.
\item[\textbf{No DIV} --] Individual diversity contributions are not counted in the biased fitness.
\item[\myblue{\textbf{No Rep}} --] \myblue{The Repair operator is not applied.}
\item[\textbf{Pop $\downarrow$} --] Smaller population: $(\mu^{elite}, \mu, \lambda)=(4,12,20)$.
\item[\textbf{Pop $\uparrow$} --] Larger population: $(\mu^{elite}, \mu, \lambda)=(16,50,80)$.
\item[\textbf{\myblue{Feas $\uparrow$}} --] \myblue{50\% feasible solutions as a target: $(\xi^\textsc{min},\xi^\textsc{max})=(0.4^{\frac{1}{1+K}},0.6^{\frac{1}{1+K}})$.}
\item[\textbf{\myblue{Feas $\uparrow\uparrow$}} --]  \myblue{75\% feasible solutions as a target: $(\xi^\textsc{min},\xi^\textsc{max})=(0.65^{\frac{1}{1+K}},0.85^{\frac{1}{1+K}})$.}\\
\end{itemize}

Each of these algorithm configurations was run 10 times on every benchmark instance for the VRP-SL, VRPPFCC and CPTP.
Table \ref{tab:sensitivity} presents the average percentage gap, best percentage gap and time of each method for each set of instances.

\begin{table}[htbp]
\vspace*{0.3cm} 
\renewcommand{\arraystretch}{1.2}
\setlength{\tabcolsep2.5pt}
\centering
\scalebox{0.83}
{
\begin{tabular}{|l|ccc|ccc|ccc|ccc|}
\hline
& \multicolumn{3}{c|}{\textbf{VRP-SL \ Set S1}} & \multicolumn{3}{c|}{\textbf{VRP-SL \  Set S2}} & \multicolumn{3}{c|}{\textbf{VRPPFCC}}& \multicolumn{3}{c|}{\textbf{CPTP}} \\
& Best-10 & Avg-10 & T(s) & Best-10 & Avg-10 &T(s) & Best-10 & Avg-10 & T(s)& Best-10 & Avg-10 & T(s) \\
\hline
\textbf{Standard}   & 0.00 & \textbf{0.01} & 16.01 & 0.00 & \textbf{0.01} & 39.66 & 0.06 & 0.21 & 340.00 & 0.00 & 0.02 & 14.07 \\ 
\textbf{EOX}   & 0.02 & 0.05 & 19.98 & 0.03 & 0.10 & 66.61 & 1.69 & 2.04 & 225.45 & 0.02 & 0.06 & 20.39 \\ 
\textbf{No SL}   & 0.00 & \textbf{0.01} & 16.70 & 0.00 & 0.02 & 41.68 & 0.07 & \textbf{0.20} & 360.61 & 0.02 & 0.04 & 15.22 \\ 
\textbf{No INF}    & 0.05 & 0.12 & 15.05 & 0.07 & 0.14 & 40.84 & 0.15 & 0.38 & 344.67 & 0.02 & 0.05 & 12.70 \\ 
\textbf{No DIV}    & 0.01 & 0.05 & 13.62 & 0.01 & 0.08 & 29.04 & 0.12 & 0.37 & 155.33 & 0.02 & 0.17 & 11.68 \\ 
\textbf{No Repair}  & 0.00 & 0.02 & 12.51 & 0.01 & 0.03 & 34.85 & 0.05 & \textbf{0.20} & 337.18 & 0.00 & 0.02 & 12.40 \\ 
\textbf{Pop $\downarrow$}  & 0.00 & \textbf{0.01} & 12.99 & 0.00 & 0.02 & 31.16 & 0.05 & 0.23 & 343.94 & 0.00 & 0.03 & 11.24 \\ 
\textbf{Pop $\uparrow$}  & 0.00 & \textbf{0.01} & 26.48 & 0.00 & 0.02 & 60.53 & 0.12 & 0.26 & 545.01 & 0.00 & \textbf{0.01} & 25.14 \\ 
\textbf{Feas $\uparrow$}  & 0.00 & \textbf{0.01} & 15.12 & 0.00 & 0.02 & 40.01 & 0.10 & 0.26 & 348.12 & 0.01 & 0.03 & 14.06 \\ 
\textbf{Feas $\uparrow\uparrow$}  & 0.01 & 0.04 & 14.93 & 0.01 & 0.05 & 40.61 & 0.16 & 0.37 & 380.05 & 0.01 & 0.05 & 14.23 \\ 
\hline
\end{tabular}
}
\caption{Sensitivity Analysis on the components of the HGS} \label{tab:sensitivity}
\end{table}

From these experiments, it appears that the proposed method is a sort of ``local optimum'' in terms of design choice and parameter settings, in the sense that any change of its main operators and parameters impacts negatively the method performance.
Still, some design choices have a much larger impact than others.
In particular, using the EOX crossover operator strongly deteriorates the method performance for VRPPFCC instances, while speeding-up the resolution. Such a speed-up may be a symptom of premature convergence due, in this case, to the crossover.
\myblue{Both diversity management and infeasible-solution management contribute significantly to the performance of the method (configurations \textbf{No INF} and \textbf{No DIV}). This confirms the earlier observations of \cite{Vidal2012,Vidal2013a}. Still, although the management of infeasible solutions is critical, deactivating the repair operator or changing the target level of feasible solutions has little impact (configurations \textbf{No Rep}, \textbf{Feas $\uparrow$} and \textbf{Feas $\uparrow\uparrow$}).
The control of the service levels in the crossover has a beneficial effect on performance for the CPTP (configuration \textbf{No SL}).
Finally, HGS is relatively insensible to reasonable changes of population size (configurations \textbf{Pop $\downarrow$} and \textbf{Pop $\uparrow$})}.

\section{Conclusions}
\label{section:concl}

In this article, we have introduced the VRP-SL, an important VRP variant arising in collaborative logistics operations, which aims to take into account the requirements of various partners via service level constraints on groups of deliveries. To establish a basis for further research, we introduced a first set of benchmark instances, a compact mathematical formulation, \myblue{a branch-and-price algorithm} and a first effective hybrid genetic search. The service level constraints tend to make the selection of services more complex, and thus new problem-tailored search operators, solution representation, crossover, LS moves, and penalty management strategies were introduced in HGS. Thanks to these elements, the proposed heuristic finds all known optimal solutions for the VRP-SL, and outperforms previous algorithms for two important special cases, the VRPPFCC and the CPTP, which have been intensively studied in past literature.
\myblue{Finally, the proposed branch-and-price algorithm was able to produce tight bounds for all problems and new optimality certificates for 63 instances, outperforming the existing exact approaches.}

The research perspectives are numerous. First, the new algorithms can still be improved via the addition of new families of cuts, better relaxations, new neighborhoods and other heuristic strategies. Moreover, the VRP-SL is only a simplification of an intricate real-life application, and the assumptions about the time constraints, the dynamics and stochasticity of the problem were voluntarily simplified to allow for reproducibility. Guaranteeing, in a stochastic and on-line context, the satisfaction of contractual service levels for prize-collecting problems is an important challenge, as the violation of such obligations can lead to large penalties or lost contracts. To circumvent this risk while mitigating the cost of robustness, it is possible to search for robust solutions on a larger planning horizon, and consider alternative transportation modes as a recourse. These aspects, and the interactions between them, will be considered in future works.


\section{Acknowledgments}
 \myblue{This research is partially supported by CNPq and CAPES in Brazil, and the National Foundation for Science and Technology Development (NAFOSTED) in Vietnam.}


\begin{thebibliography}{69}
\expandafter\ifx\csname natexlab\endcsname\relax\def\natexlab#1{#1}\fi
\expandafter\ifx\csname url\endcsname\relax
  \def\url#1{{\tt #1}}\fi
\expandafter\ifx\csname urlprefix\endcsname\relax\def\urlprefix{URL }\fi
\expandafter\ifx\csname urlstyle\endcsname\relax
  \expandafter\ifx\csname doi\endcsname\relax
  \def\doi#1{doi:\discretionary{}{}{}#1}\fi \else
  \expandafter\ifx\csname doi\endcsname\relax
  \def\doi{doi:\discretionary{}{}{}\begingroup \urlstyle{rm}\Url}\fi \fi

\bibitem[{Archetti et~al.(2013)Archetti, Bianchessi, and
  Speranza}]{Archetti2013a}
Archetti, C., N.~Bianchessi, M.G. Speranza. 2013.
\newblock {Optimal solutions for routing problems with profits}.
\newblock {\it Discrete Applied Mathematics\/} {\bf 161}(4-5) 547--557.

\bibitem[{Archetti et~al.(2009)Archetti, Feillet, Hertz, and
  Speranza}]{Archetti2008d}
Archetti, C., D.~Feillet, A.~Hertz, M.G. Speranza. 2009.
\newblock {The capacitated team orienteering and profitable tour problems}.
\newblock {\it Journal of the Operational Research Society\/} {\bf 60}(6)
  831--842.

\bibitem[{Archetti et~al.(2014)Archetti, Speranza, and Vigo}]{Archetti2014}
Archetti, C., M.G. Speranza, D.~Vigo. 2014.
\newblock {Vehicle routing problems with profits}.
\newblock P.~Toth, D.~Vigo, eds., {\it Vehicle Routing: Problems, Methods, and
  Applications\/}. SIAM, Philadelphia, PA, 273--297.

\bibitem[{Augerat et~al.(1995)Augerat, Belenguer, Benavent, Coberan, Naddef,
  and Rinaldi}]{Augerat1995}
Augerat, P., J.M. Belenguer, E.~Benavent, A.~Coberan, D.~Naddef, G.~Rinaldi.
  1995.
\newblock {Computational results with a branch-and-cut code for the capacitated
  vehicle routing problem}.
\newblock Tech. rep., Universit{\'{e}} Joseph Fourier, Grenoble, France.

\bibitem[{Baldacci et~al.(2009)Baldacci, Bartolini, and
  Laporte}]{Baldacci2009a}
Baldacci, R., E.~Bartolini, G.~Laporte. 2009.
\newblock {Some applications of the generalized vehicle routing problem}.
\newblock {\it Journal of the Operational Research Society\/} {\bf 61}(7)
  1072--1077.

\bibitem[{Baldacci et~al.(2008)Baldacci, Christofides, and
  Mingozzi}]{baldacci-2008}
Baldacci, R., N.~Christofides, A.~Mingozzi. 2008.
\newblock {An exact algorithm for the vehicle routing problem based on the set
  partitioning formulation with additional cuts}.
\newblock {\it Mathematical Programming\/} {\bf 115}(2) 351--385.

\bibitem[{Baldacci et~al.(2004)Baldacci, Hadjiconstantinou, and
  Mingozzi}]{Baldacci2004}
Baldacci, R., E.~Hadjiconstantinou, A.~Mingozzi. 2004.
\newblock {An exact algorithm for the capacitated vehicle routing problem based
  on a two-commodity network flow formulation}.
\newblock {\it Operations Research\/} {\bf 52}(5) 723--738.

\bibitem[{Baldacci et~al.(2011)Baldacci, Mingozzi, and Roberti}]{baldacci-2011}
Baldacci, R., A.~Mingozzi, R.~Roberti. 2011.
\newblock New route relaxation and pricing strategies for the vehicle routing
  problem.
\newblock {\it Operations Research\/} {\bf 59}(5) 1269--1283.

\bibitem[{Beasley(1983)}]{Beasley1983}
Beasley, J.E. 1983.
\newblock {Route first-cluster second methods for vehicle routing}.
\newblock {\it Omega\/} {\bf 11}(4) 403--408.

\bibitem[{Bektas et~al.(2011)Bektas, Erdogan, and Ropke}]{Bektas2011}
Bektas, T., G.~Erdogan, S.~Ropke. 2011.
\newblock {Formulations and branch-and-cut algorithms for the generalized
  vehicle routing problem}.
\newblock {\it Transportation Science\/} {\bf 45}(3) 299--316.

\bibitem[{B{\'{e}}rub{\'{e}} et~al.(2009)B{\'{e}}rub{\'{e}}, Gendreau, and
  Potvin}]{Berube2009}
B{\'{e}}rub{\'{e}}, J.F., M.~Gendreau, J.Y. Potvin. 2009.
\newblock {A branch-and-cut algorithm for the undirected prize collecting
  traveling salesman problem}.
\newblock {\it Networks\/} {\bf 54}(1) 56--67.

\bibitem[{Bolduc et~al.(2008)Bolduc, Renaud, Boctor, and Laporte}]{Bolduc2008}
Bolduc, M.-C., J.~Renaud, F.~Boctor, G.~Laporte. 2008.
\newblock {A perturbation metaheuristic for the vehicle routing problem with
  private fleet and common carriers}.
\newblock {\it Journal of the Operational Research Society\/} {\bf 59}(6)
  776--787.

\bibitem[{Boussier et~al.(2007)Boussier, Feillet, and Gendreau}]{Boussier2007}
Boussier, S., D.~Feillet, M.~Gendreau. 2007.
\newblock {An exact algorithm for team orienteering problems}.
\newblock {\it 4OR\/} {\bf 5}(3) 211--230.

\bibitem[{Campos et~al.(2014)Campos, Mart{\'{i}}, S{\'{a}}nchez-Oro, and
  Duarte}]{Campos2014}
Campos, V., R.~Mart{\'{i}}, J.~S{\'{a}}nchez-Oro, A.~Duarte. 2014.
\newblock {GRASP with path relinking for the orienteering problem}.
\newblock {\it Journal of the Operational Research Society\/} {\bf 65}
  1800--1813.

\bibitem[{Chaves and Lorena(2008)}]{Chaves2008}
Chaves, A.A., L.A.N. Lorena. 2008.
\newblock {Hybrid metaheuristic for the prize collecting travelling salesman
  problem}.
\newblock J.~van Hemert, C.~Cotta, eds., {\it Evolutionary Computation in
  Combinatorial Optimization\/}, {\it LNCS\/}, vol. 4972. Springer, Berlin,
  Heidelberg, 123--134.

\bibitem[{Christofides et~al.(1981)Christofides, Mingozzi, and
  Toth}]{christofides-1981}
Christofides, N., A.~Mingozzi, P.~Toth. 1981.
\newblock {Exact algorithms for the vehicle routing problem, based on spanning
  tree and shortest path relaxations}.
\newblock {\it Mathematical Programming\/} {\bf 20} 255--282.

\bibitem[{Contardo and Martinelli(2014)}]{contardo-2014}
Contardo, C., R.~Martinelli. 2014.
\newblock A new exact algorithm for the multi-depot vehicle routing problem
  under capacity and route length constraints.
\newblock {\it Discrete Optimization\/} {\bf 12} 129--146.

\bibitem[{C{\^{o}}t{\'{e}} and Potvin(2009)}]{Cote2009a}
C{\^{o}}t{\'{e}}, J.-F., J.-Y. Potvin. 2009.
\newblock {A tabu search heuristic for the vehicle routing problem with private
  fleet and common carrier}.
\newblock {\it European Journal of Operational Research\/} {\bf 198}(2)
  464--469.

\bibitem[{Dang et~al.(2013)Dang, Guibadj, and Moukrim}]{Dang2013}
Dang, D.-C., R.N. Guibadj, A.~Moukrim. 2013.
\newblock {An effective PSO-inspired algorithm for the team orienteering
  problem}.
\newblock {\it European Journal of Operational Research\/} {\bf 229}(2)
  332--344.

\bibitem[{Dell'Amico et~al.(1995)Dell'Amico, Maffioli, and
  V{\"{a}}rbrand}]{DellAmico1995}
Dell'Amico, M., F.~Maffioli, P.~V{\"{a}}rbrand. 1995.
\newblock {On prize-collecting tours and the asymmetric travelling salesman
  problem}.
\newblock {\it International Transactions in Operational Research\/} {\bf 2}(3)
  297--308.

\bibitem[{Dror(1994)}]{dror-1994}
Dror, M. 1994.
\newblock {Note on the complexity of the shortest path models for column
  generation in VRPTW}.
\newblock {\it Operations Research\/} {\bf 42}(5) 977--978.

\bibitem[{El-Hajj et~al.(2016)El-Hajj, Dang, and Moukrim}]{ElHajj2016}
El-Hajj, R., D.-C. Dang, A.~Moukrim. 2016.
\newblock {Solving the team orienteering problem with cutting planes}.
\newblock {\it Computers {\&} Operations Research\/} {\bf 74} 21--30.

\bibitem[{Feillet et~al.(2005)Feillet, Dejax, and Gendreau}]{Feillet2005}
Feillet, D., P.~Dejax, M.~Gendreau. 2005.
\newblock {Traveling salesman problems with profits}.
\newblock {\it Transportation Science\/} {\bf 39}(2) 188--205.

\bibitem[{Fischetti et~al.(1998)Fischetti, Gonzalez, and Toth}]{Fischetti1998}
Fischetti, M., J.J.S. Gonzalez, P.~Toth. 1998.
\newblock {Solving the orienteering problem through branch-and-cut}.
\newblock {\it INFORMS Journal on Computing\/} {\bf 10}(2) 133--148.

\bibitem[{Fischetti and Toth(1988)}]{Fischetti1988}
Fischetti, M., P.~Toth. 1988.
\newblock {An additive approach for the optimal solution of the
  prize-collecting travelling salesman problem}.
\newblock B.L. Golden, A.A. Assad, eds., {\it Vehicle Routing: Methods and
  Studies\/}. Elsevier, 319--343.

\bibitem[{Fukasawa et~al.(2006)Fukasawa, Longo, Lysgaard, {Poggi de
  Arag{{\~a}}o}, Reis, Uchoa, and Werneck}]{fukasawa-2006}
Fukasawa, R., H.~Longo, J.~Lysgaard, M.~{Poggi de Arag{{\~a}}o}, M.~Reis,
  E.~Uchoa, R.F. Werneck. 2006.
\newblock {Robust branch-and-cut-and-price for the capacitated vehicle routing
  problem}.
\newblock {\it Mathematical Programming\/} {\bf 106}(3) 491--511.

\bibitem[{Gavish and Graves(1978)}]{Gavish1978}
Gavish, B., S.C. Graves. 1978.
\newblock {The travelling salesman problem and related problems}.
\newblock Tech. rep., Massachusetts Institute of Technology, Cambridge, MA,
  USA.

\bibitem[{Gendreau et~al.(1998)Gendreau, Laporte, and Semet}]{Gendreau1998b}
Gendreau, M., G.~Laporte, F.~Semet. 1998.
\newblock {A tabu search heuristic for the undirected selective travelling
  salesman problem}.
\newblock {\it European Journal of Operational Research\/} {\bf 106}(2-3)
  539--545.

\bibitem[{Ghiani and Improta(2000)}]{Ghiani2000}
Ghiani, G., G.~Improta. 2000.
\newblock {An efficient transformation of the generalized vehicle routing
  problem}.
\newblock {\it European Journal of Operational Research\/} {\bf 122}(1) 11--17.

\bibitem[{H{\`{a}} et~al.(2013)H{\`{a}}, Bostel, Langevin, and
  Rousseau}]{Ha2013}
H{\`{a}}, M.H., N.~Bostel, A.~Langevin, L.-M. Rousseau. 2013.
\newblock {An exact algorithm and a metaheuristic for the multi-vehicle
  covering tour problem with a constraint on the number of vertices}.
\newblock {\it European Journal of Operational Research\/} {\bf 226}(2)
  211--220.

\bibitem[{H{\`{a}} et~al.(2014)H{\`{a}}, Bostel, Langevin, and
  Rousseau}]{Ha2014}
H{\`{a}}, M.H., N.~Bostel, A.~Langevin, L.-M. Rousseau. 2014.
\newblock {An exact algorithm and a metaheuristic for the generalized vehicle
  routing problem with flexible fleet size}.
\newblock {\it Computers {\&} Operations Research\/} {\bf 43}(1) 9--19.

\bibitem[{Huijink et~al.(2014)Huijink, Goos, and Peeters}]{Huijink2014}
Huijink, S., K.~Goos, R.~Peeters. 2014.
\newblock {An adaptable variable neighborhood search for the vehicle routing
  problem with order outsourcing}.
\newblock Tech. rep., Tilburg University.

\bibitem[{Irnich and Villeneuve(2006)}]{irnich-2006}
Irnich, S., D.~Villeneuve. 2006.
\newblock {The shortest-path problem with resource constraints and k-cycle
  elimination for k $\geq$ 3}.
\newblock {\it INFORMS Journal on Computing\/} {\bf 18}(3) 391--406.

\bibitem[{Jayanth and Keah-Choon(2010)}]{jayaram-2010}
Jayanth, J., T.~Keah-Choon. 2010.
\newblock Supply chain integration with third-party logistics providers.
\newblock {\it International Journal of Production Economics\/} {\bf 125}(2)
  262--271.

\bibitem[{Jepsen et~al.(2014)Jepsen, Petersen, Spoorendonk, and
  Pisinger}]{Jepsen2014}
Jepsen, M.K., B.~Petersen, S.~Spoorendonk, D.~Pisinger. 2014.
\newblock {A branch-and-cut algorithm for the capacitated profitable tour
  problem}.
\newblock {\it Discrete Optimization\/} {\bf 14}(1) 78--96.

\bibitem[{Kataoka and Morito(1988)}]{Kataoka1988}
Kataoka, S., S.~Morito. 1988.
\newblock {An algorithm for single constraint maximum collection problem}.
\newblock {\it Journal of the Operations Research Society of Japan\/} {\bf
  31}(4) 515--531.

\bibitem[{Ke et~al.(2008)Ke, Archetti, and Feng}]{Ke2008}
Ke, L., C.~Archetti, Z.~Feng. 2008.
\newblock {Ants can solve the team orienteering problem}.
\newblock {\it Computers {\&} Industrial Engineering\/} {\bf 54}(3) 648--665.

\bibitem[{Ke et~al.(2015)Ke, Zhai, Li, and Chan}]{Ke2015}
Ke, L., L.~Zhai, J.~Li, F.T.S. Chan. 2015.
\newblock {Pareto mimic algorithm: An approach to the team orienteering
  problem}.
\newblock {\it Omega\/} {\bf 61}(1) 155--166.

\bibitem[{Keshtkaran et~al.(2016)Keshtkaran, Ziarati, Bettinelli, and
  Vigo}]{Keshtkaran2015}
Keshtkaran, M., K.~Ziarati, A.~Bettinelli, D.~Vigo. 2016.
\newblock {Enhanced exact solution methods for the Team Orienteering Problem}.
\newblock {\it International Journal of Production Research\/} {\bf 54}(2)
  591--601.

\bibitem[{Kim et~al.(2013)Kim, Li, and Johnson}]{Kim2013}
Kim, B.-I., H.~Li, A.L. Johnson. 2013.
\newblock {An augmented large neighborhood search method for solving the team
  orienteering problem}.
\newblock {\it Expert Systems with Applications\/} {\bf 40}(8) 3065--3072.

\bibitem[{Letchford and Salazar-Gonz{\'{a}}lez(2006)}]{Letchford2006a}
Letchford, A.N., J.-J. Salazar-Gonz{\'{a}}lez. 2006.
\newblock {Projection results for vehicle routing}.
\newblock {\it Mathematical Programming\/} {\bf 105}(2-3) 251--274.

\bibitem[{Leuschner et~al.(2014)Leuschner, Carter, Goldsby, and
  Rogers}]{leuschner-2014}
Leuschner, R., C.R. Carter, T.J. Goldsby, Z.S. Rogers. 2014.
\newblock Third-party logistics: A meta-analytic review and investigation of
  its impact on performance.
\newblock {\it Journal of Supply Chain Management\/} {\bf 50}(1) 21--43.

\bibitem[{Li and Tian(2016)}]{Li2016}
Li, K., H.~Tian. 2016.
\newblock {A two-level self-adaptive variable neighborhood search algorithm for
  the prize-collecting vehicle routing problem}.
\newblock {\it Applied Soft Computing\/} {\bf 43}(1) 469--479.

\bibitem[{Lin(2013)}]{Lin2013}
Lin, S.-W. 2013.
\newblock {Solving the team orienteering problem using effective multi-start
  simulated annealing}.
\newblock {\it Applied Soft Computing\/} {\bf 13}(2) 1064--1073.

\bibitem[{Lopez et~al.(1998)Lopez, Carter, and Gendreau}]{Lopez1998}
Lopez, L., M.W. Carter, M.~Gendreau. 1998.
\newblock {The hot strip mill production scheduling problem: A tabu search
  approach}.
\newblock {\it European Journal of Operational Research\/} {\bf 106}(2-3)
  317--335.

\bibitem[{Martinelli et~al.(2014)Martinelli, Pecin, and
  Poggi}]{martinelli-2014}
Martinelli, R., D.~Pecin, M.~Poggi. 2014.
\newblock Efficient elementary and restricted non-elementary route pricing.
\newblock {\it European Journal of Operational Research\/} {\bf 239}(1)
  102--111.

\bibitem[{Martinelli et~al.(2011)Martinelli, Pecin, Poggi, and
  Longo}]{martinelli-2011}
Martinelli, R., D.~Pecin, M.~Poggi, H.~Longo. 2011.
\newblock A branch-cut-and-price algorithm for the capacitated arc routing
  problem.
\newblock P.~Pardalos, S.~Rebennack, eds., {\it Experimental Algorithms\/},
  {\it LNCS\/}, vol. 6630. Springer, 315--326.

\bibitem[{M{\"{u}}hlenbein and Schlierkamp-Voosen(1993)}]{Muhlenbein1993}
M{\"{u}}hlenbein, H., D.~Schlierkamp-Voosen. 1993.
\newblock {Predictive models for the breeder genetic algorithm I. Continuous
  parameter optimization}.
\newblock {\it Evolutionary Computation\/} {\bf 1}(1) 25--49.

\bibitem[{Pedro et~al.(2013)Pedro, Saldanha, and Camargo}]{Pedro2013}
Pedro, O., R.~Saldanha, R.~Camargo. 2013.
\newblock {A tabu search approach for the prize collecting traveling salesman
  problem}.
\newblock {\it Electronic Notes in Discrete Mathematics\/} {\bf 41}(1)
  261--268.

\bibitem[{Pessoa et~al.(2010)Pessoa, Uchoa, {Poggi de Arag{{\~a}}o}, and
  Rodrigues}]{pessoa-2010}
Pessoa, A., E.~Uchoa, M.~{Poggi de Arag{{\~a}}o}, R.~Rodrigues. 2010.
\newblock Exact algorithm over an arc-time-indexed formulation for parallel
  machine scheduling problems.
\newblock {\it Mathematical Programming Computation\/} {\bf 2}(3) 259--290.

\bibitem[{Potvin and Naud(2011)}]{Potvin2011}
Potvin, J.-Y., M.-A. Naud. 2011.
\newblock {Tabu search with ejection chains for the vehicle routing problem
  with private fleet and common carrier}.
\newblock {\it Journal of the Operational Research Society\/} {\bf 62}(2)
  326--336.

\bibitem[{Prins(2004)}]{Prins2004}
Prins, C. 2004.
\newblock {A simple and effective evolutionary algorithm for the vehicle
  routing problem}.
\newblock {\it Computers {\&} Operations Research\/} {\bf 31}(12) 1985--2002.

\bibitem[{Schilde et~al.(2009)Schilde, Doerner, Hartl, and
  Kiechle}]{Schilde2009}
Schilde, M., K.F. Doerner, R.F. Hartl, G.~Kiechle. 2009.
\newblock {Metaheuristics for the bi-objective orienteering problem}.
\newblock {\it Swarm Intelligence\/} {\bf 3}(1) 179--201.

\bibitem[{Souffriau et~al.(2010)Souffriau, Vansteenwegen, {Vanden Berghe}, and
  {Van Oudheusden}}]{Souffriau2010}
Souffriau, W., P.~Vansteenwegen, G.~{Vanden Berghe}, D.~{Van Oudheusden}. 2010.
\newblock {A path relinking approach for the team orienteering problem}.
\newblock {\it Computers {\&} Operations Research\/} {\bf 37}(11) 1853--1859.

\bibitem[{Stefansson(2006)}]{stefansson-2006}
Stefansson, G. 2006.
\newblock Collaborative logistics management and the role of third-party
  service providers.
\newblock {\it International Journal of Physical Distribution \& Logistics
  Management\/} {\bf 36}(2) 76--92.

\bibitem[{Stenger et~al.(2013)Stenger, Schneider, and Goeke}]{Stenger2013}
Stenger, A., M.~Schneider, D.~Goeke. 2013.
\newblock {The prize-collecting vehicle routing problem with single and
  multiple depots and non-linear cost}.
\newblock {\it EURO Journal on Transportation and Logistics\/} {\bf 2}(1-2)
  57--87.

\bibitem[{Stenger et~al.(2012)Stenger, Vigo, Enz, and Schwind}]{Stenger2012a}
Stenger, A., D.~Vigo, S.~Enz, M.~Schwind. 2012.
\newblock {An adaptive variable neighborhood search algorithm for a vehicle
  routing problem arising in small package shipping}.
\newblock {\it Transportation Science\/} {\bf 47}(1) 64--80.

\bibitem[{Tang and Wang(2006)}]{Tang2006a}
Tang, L., X.~Wang. 2006.
\newblock {Iterated local search algorithm based on very large-scale
  neighborhood for prize-collecting vehicle routing problem}.
\newblock {\it The International Journal of Advanced Manufacturing
  Technology\/} {\bf 29}(11-12) 1246--1258.

\bibitem[{Vansteenwegen et~al.(2009)Vansteenwegen, Souffriau, Berghe, and
  Oudheusden}]{Vansteenwegen2009b}
Vansteenwegen, P., W.~Souffriau, G.V. Berghe, D.V. Oudheusden. 2009.
\newblock {A guided local search metaheuristic for the team orienteering
  problem}.
\newblock {\it European Journal of Operational Research\/} {\bf 196}(1)
  118--127.

\bibitem[{Vansteenwegen et~al.(2010)Vansteenwegen, Souffriau, and
  Oudheusden}]{Vansteenwegen2010}
Vansteenwegen, P., W.~Souffriau, D.V. Oudheusden. 2010.
\newblock {The orienteering problem: A survey}.
\newblock {\it European Journal of Operational Research\/} {\bf 209}(1) 1--10.

\bibitem[{Vidal(2016)}]{Vidal2016}
Vidal, T. 2016.
\newblock {Technical note: Split algorithm in O(n) for the capacitated vehicle
  routing problem}.
\newblock {\it Computers {\&} Operations Research\/} {\bf 69} 40--47.

\bibitem[{Vidal et~al.(2012)Vidal, Crainic, Gendreau, Lahrichi, and
  Rei}]{Vidal2012}
Vidal, T., T.G. Crainic, M.~Gendreau, N.~Lahrichi, W.~Rei. 2012.
\newblock {A hybrid genetic algorithm for multidepot and periodic vehicle
  routing problems}.
\newblock {\it Operations Research\/} {\bf 60}(3) 611--624.

\bibitem[{Vidal et~al.(2013)Vidal, Crainic, Gendreau, and Prins}]{Vidal2012a}
Vidal, T., T.G. Crainic, M.~Gendreau, C.~Prins. 2013.
\newblock {Heuristics for multi-attribute vehicle routing problems: A survey
  and synthesis}.
\newblock {\it European Journal of Operational Research\/} {\bf 231}(1) 1--21.

\bibitem[{Vidal et~al.(2014)Vidal, Crainic, Gendreau, and Prins}]{Vidal2012b}
Vidal, T., T.G. Crainic, M.~Gendreau, C.~Prins. 2014.
\newblock {A unified solution framework for multi-attribute vehicle routing
  problems}.
\newblock {\it European Journal of Operational Research\/} {\bf 234}(3)
  658--673.

\bibitem[{Vidal et~al.(2015)Vidal, Crainic, Gendreau, and Prins}]{Vidal2013a}
Vidal, T., T.G. Crainic, M.~Gendreau, C.~Prins. 2015.
\newblock {Time-window relaxations in vehicle routing heuristics}.
\newblock {\it Journal of Heuristics\/} {\bf 21}(3) 329--358.

\bibitem[{Vidal et~al.(2016)Vidal, Maculan, Ochi, and Penna}]{Vidal2014}
Vidal, T., N.~Maculan, L.S. Ochi, P.H.V. Penna. 2016.
\newblock {Large neighborhoods with implicit customer selection for vehicle
  routing problems with profits}.
\newblock {\it Transportation Science\/} {\bf 50}(2) 720--734.

\bibitem[{Whitley et~al.(1989)Whitley, Starkweather, and Fuquay}]{Whitley1989}
Whitley, L.D., T.~Starkweather, D.~Fuquay. 1989.
\newblock {Scheduling problems and traveling salesmen: The genetic edge
  recombination operator}.
\newblock {\it Proceedings of the 3rd International Conference on Genetic
  Algorithms\/}. Morgan Kaufmann, San Francisco, CA, USA, 133--140.

\bibitem[{Yadollahpour et~al.(2009)Yadollahpour, Bijari, Kavosh, and
  Mahnam}]{Yadollahpour2009}
Yadollahpour, M.R., M.~Bijari, S.~Kavosh, M.~Mahnam. 2009.
\newblock {Guided local search algorithm for hot strip mill scheduling problem
  with considering hot charge rolling}.
\newblock {\it International Journal of Advanced Manufacturing Technology\/}
  {\bf 45}(11) 1215--1231.

\bibitem[{Zhang et~al.(2009)Zhang, Chaovalitwongse, Zhang, and
  Pardalos}]{Zhang2009}
Zhang, T., W.~Chaovalitwongse, Y.-J. Zhang, P.M. Pardalos. 2009.
\newblock {The hot-rolling batch scheduling method based on the prize
  collecting vehicle routing problem}.
\newblock {\it Journal of Industrial and Management Optimization\/} {\bf 5}(4)
  749--765.

\end{thebibliography}

\end{document}